\documentclass[11pt,letter]{article}

\usepackage{amsmath,amsfonts,amsthm,amssymb}
\usepackage{hyperref}
\usepackage[english]{babel}
\usepackage{mathrsfs}
\usepackage{color}
\usepackage[left=2.5cm, top=2.5cm, bottom=2.5cm, right=2.5cm]{geometry}

\newtheorem{theorem}{Theorem}
\newtheorem{assumption}[theorem]{Assumption}

\newtheorem{corollary}[theorem]{Corollary}

\newtheorem{lemma}[theorem]{Lemma}
\newtheorem{proposition}[theorem]{Proposition}

\theoremstyle{definition}

\newtheorem{example}[theorem]{Example}
\newtheorem{remark}[theorem]{Remark}

\numberwithin{equation}{section}
\numberwithin{theorem}{section}

\definecolor{Red}{rgb}{1.00, 0.00, 0.00}

\definecolor{Green}{rgb}{0.2, 0.6, 0.2}

\definecolor{Blue}{rgb}{0.00, 0.00, 1.00}

\definecolor{DRed}{rgb}{0.7, 0.3, 0.00}

\def\convdist{\;{\stackrel{\mathfrak{D}}{\longrightarrow}}\;}

\title{Lower Bounds on the Generalized Central Moments of the Optimal Alignments Score of Random Sequences}

\author{Ruoting Gong\thanks{Department of Applied Mathematics, Illinois Institute of Technology, Chicago, IL 60616, U.S.A. {\tt Email:\,rgong2@iit.edu}.}
\and Christian Houdr\'{e}\thanks{School of Mathematics, Georgia Institute of Technology, Atlanta, GA 30332, U.S.A. {\tt
Email:\,houdre@math.gatech.edu}. Research supported in part by the grant \# 246283 from the Simons Foundation and by a Simons Fellowship grant \# 267336.}
\and J\"{u}ri Lember\thanks{Institute of Mathematical Statistics, University of Tartu, Liivi 2-513 50409, Tartu, Estonia. {\tt Email:\,juril@ut.ee}. Research supported by Estonian Science foundation Grant no. 5822 and by institutional research funding IUT34-5.}}

\date{November 24, 2016}

\begin{document}

\maketitle

\begin{abstract}

We present a general approach to the problem of determining tight asymptotic lower bounds for generalized central moments of the optimal alignment score of two independent sequences of i.i.d.~random variables. At first, these are obtained under a main assumption for which sufficient conditions are provided. When the main assumption fails, we nevertheless develop a ``uniform approximation" method leading to asymptotic lower bounds. Our general results are then applied to the length of the longest common subsequence of binary strings, in which case asymptotic lower bounds are obtained for the moments and the exponential moments of the optimal score. As a byproduct, a local upper bound on the rate function associated with the length of the longest common subsequences of two binary strings is also obtained.

\vspace{0.3 cm}


\vspace{0.3 cm}


\end{abstract}

\section{Introduction}

Throughout this paper, $\mathbf{X}_{n}:=(X_{1},X_{2},\ldots,X_{n})$ and $\mathbf{Y}_{n}:=(Y_{1},Y_{2},\ldots,Y_{n})$ are two random strings, usually referred to as the (finite) sequences, so that any random variable $X_{i}$ or $Y_{i}$, $i=1,\ldots,n$, takes its values in a fixed finite alphabet $\mathcal{A}$. The sequences $\mathbf{X}_{n}$ and $\mathbf{Y}_{n}$ are assumed to have the same distribution and to be independent. The sample space of $\mathbf{X}_{n}$ and $\mathbf{Y}_{n}$ will be denoted by $\mathcal{X}_{n}$; clearly $\mathcal{X}_{n}\subseteq\mathcal{A}^{n}$, but, depending on the model, equality might not hold.

The problem of measuring the similarity of $\mathbf{X}_{n}$ and $\mathbf{Y}_{n}$ is central to many areas of applications including computational molecular biology (cf.~\cite{ChristianiniHahn:2007}, \cite{DurbinEddyKroghMitchison:1998}, \cite{Pevzner:2000}, \cite{SmithWaterman:1981} and~\cite{Waterman:1995}) and computational linguistics (cf.~\cite{YangLi:2003}, \cite{LinOch:2004}, \cite{Melamed:1995} and~\cite{Melamed:1999}). In this paper we adopt the notation of~\cite{LemberMatzingerTorres:2012}, namely, we consider a general scoring scheme, where $S:\mathcal{A}\times\mathcal{A}\rightarrow\mathbb{R}^{+}$ is a symmetric {\it pairwise scoring function}, that assigns a score to each pair of letters from $\mathcal{A}$. Throughout, an {\it alignment} is a pair $(\boldsymbol\pi,\boldsymbol\mu)$, where $\boldsymbol\pi:=(\pi_{1},\pi_{2},\ldots,\pi_{k})$ and $\boldsymbol\mu:=(\mu_{1},\mu_{2},\ldots,\mu_{k})$ are two increasing sequences of positive integers such that $1\leq\pi_{1}<\pi_{2}<\cdots<\pi_{k}\leq n$ and $1\leq\mu_{1}<\mu_{2}<\cdots<\mu_{k}\leq n$. The integer $k$ is the number of aligned letters, and $n-k$ is the number of {\it gaps} in the alignment. Note that our definition of gap differs from the one that is commonly used in the sequence alignment literature. More precisely, in most of the literature, a gap is a block of consecutive {\it indels} (insertion and deletion, formally denoted by one or more consecutive ``--") in both strings, and it depends on the alignment. For example, the left alignment below has four gaps, while the right one has five gaps.
\begin{align*}
\begin{array}{ccccccccc} 1 & 1 & 3 & 1 & 2 & - & - & - & 3 \\ 1 & - & 3 & - & 2 & 1 & 1 & 1 & - \end{array}\qquad\qquad\begin{array}{ccccccccc} 1 & 1 & 3 & 1 & 2 & - & 3 & - & - \\ 1 & - & 3 & - & 2 & 1 & - & 1 & 1 \end{array}\quad .
\end{align*}
Since we consider sequences of equal length, and since we do not have a gap opening penalty (which refers to a constant cost to open a gap of any length), our gap corresponds to a pair of indels, one on the $\mathbf{X}$-side and the other on the $\mathbf{Y}$-side. In other words, the number of gaps in this sense is the number of indels in either one of the sequences. Hence, in our framework, both the left and the right alignment above have three gaps.

Given a symmetric pairwise scoring function $S:\mathcal{A}\times\mathcal{A}\rightarrow\mathbb{R}^{+}$ and a {\it constant} gap price $\delta\in\mathbb{R}$, the score of the alignment $(\boldsymbol\pi,\boldsymbol\mu)$, when aligning $\mathbf{X}_{n}$ and $\mathbf{Y}_{n}$, is defined by
\begin{align*}
U_{(\boldsymbol\pi,\boldsymbol\mu)}\!\left(\mathbf{X}_{n};\mathbf{Y}_{n}\right):=\sum_{i=1}^{k}S(X_{\pi_{i}},Y_{\mu_{i}})+\delta(n-k).
\end{align*}
In our general scoring scheme $\delta$ can be positive, although usually $\delta\leq 0$ penalizes matches with ``--". For negative $\delta$, the quantity $-\delta$ is usually called the {\it gap penalty}. The optimal alignment score of $\mathbf{X}_{n}$ and $\mathbf{Y}_{n}$ is then defined as
\begin{align*}
L_{n}=L_{n}\!\left(\mathbf{X}_{n};\mathbf{Y}_{n}\right):=\max_{(\boldsymbol\pi,\boldsymbol\mu)}U_{(\boldsymbol\pi,\boldsymbol\mu)}\!\left(\mathbf{X}_{n};\mathbf{Y}_{n}\right),
\end{align*}
where the maximum is taken over all possible alignments. To simplify notations, in what follows, we often set $\mathbf{Z}_{n}:=(\mathbf{X}_{n},\mathbf{Y}_{n})$ so that $L_{n}=L_{n}(\mathbf{Z}_{n})$. When $\delta=0$ and the scoring function assigns the value one to every pair of common letters and zero to all other pairs, i.e.,
\begin{align*}
S(a,b)=\left\{\begin{array}{ll} 1, & \hbox{if $a=b$;} \\ 0, & \hbox{if $a\ne b$,} \end{array}\right.
\end{align*}
then $L_{n}(\mathbf{Z}_{n})$ is just the maximal number of pairs of aligned letters -- the length of the {\it longest common subsequences} (LCS) of $\mathbf{X}_{n}$ and $\mathbf{Y}_{n}$, which is probably the most frequently used measure of global similarity.

\subsection{A Summary of Known Results}

The study of $L_{n}$ is not only of theoretical interest but also has some practical consequences. For example, it is useful for distinguishing related (homologous) pairs of strings from unrelated ones. Unfortunately, for fixed $n$ and an arbitrary distribution for $\mathbf{X}_{n}$ and $\mathbf{Y}_{n}$, the distribution of $L_{n}$ is unknown. In view of these, an alternative approach is to study the typical values and the fluctuations of $L_{n}$. When $\mathbf{X}_{n}$ and $\mathbf{Y}_{n}$ are taken from an ergodic process, by Kingman's subadditive ergodic theorem, there exists a constant $\gamma^{*}$ such that
\begin{align}\label{eq:ConvLn}
\frac{L_{n}}{n}\rightarrow\gamma^{*}\quad\text{a.}\,\text{s.}\,\,\,\text{and in }\,L_{1},\quad\text{as }\,n\rightarrow\infty.
\end{align}
In the LCS case, the existence of $\gamma^{*}$ was first shown by Chv\'{a}tal and Sankoff~\cite{ChvatalSankoff:1975}, but its exact value (or an identity for it) remains unknown even for independent i.i.d. Bernoulli sequences (although numerically estimated). Alexander~\cite{Alexander:1994} obtained the rate of convergence in \eqref{eq:ConvLn} in the LCS case, which was extended to general scoring functions in~\cite{LemberMatzingerTorres:2012}, and to multiple (two or more) independent i.i.d. sequences with general score function in~\cite{GongHoudreIslak:2016}.

Since the exact distribution of $L_{n}$ is rather hard to determine even for moderate $n$, it is natural to look for a limiting theorem, e.g.,
\begin{align*}
\frac{L_{n}-\mathbb{E}(L_{n})}{n^{\alpha}}\convdist\mathcal{P},\quad n\rightarrow\infty,
\end{align*}
for some $\alpha\in(0,1)$. Here $\mathcal{P}$ is a limiting distribution, and $\convdist$ stands for convergence in law. Typically, one expects $\alpha=1/2$, and $\mathcal{P}$ to be a centered normal distribution, i.e.,
\begin{align}\label{eq:CLTLn}
\frac{L_{n}-\mathbb{E}(L_{n})}{\sqrt{n}}\convdist\mathcal{N},\quad n\rightarrow\infty,
\end{align}
where $\mathcal{N}$ stands for a centered normal distribution with variance $\sigma^{2}>0$. Under \eqref{eq:CLTLn}, for any $r>0$, the $r$-th absolute moment of $L_{n}$ would then be expected to grow at speed $n^{r/2}$, as $n\rightarrow\infty$. In particular, the variance would grow linearly, i.e.,
\begin{align*}
\lim_{n\rightarrow\infty}\frac{\text{Var}(L_{n})}{n}=\sigma^{2}>0,
\end{align*}
and then, \eqref{eq:CLTLn} would be equivalent to
\begin{align}\label{eq:CLTLn2}
\frac{L_{n}-\mathbb{E}(L_{n})}{\sqrt{\text{Var}(L_{n})}}\convdist\mathcal{N}(0,1),\quad n\rightarrow\infty.
\end{align}
Such a limiting theorem would allow to construct asymptotic tests, based on the optimal score $L_{n}$, hypothesizing that two given sequences $\mathbf{X}_{n}$ and $\mathbf{Y}_{n}$ are independent (not homologous) or not.

Note that in the gapless case, i.e., when $\delta=-\infty$, the optimal score is just the sum of the pairwise scores $L_{n}=\sum_{i=1}^{n}S(X_{i},Y_{i})$, and thus, under rather general assumptions on $\mathbf{X}_{n}$ and $\mathbf{Y}_{n}$, the limiting theorem \eqref{eq:CLTLn2} (equivalently, \eqref{eq:CLTLn}) holds true. This observation suggests to conjecture that \eqref{eq:CLTLn2} also holds in more general cases, in particular, in the LCS case. No such type of limiting theorem was known until~\cite{HoudreIslak:2015}, which proved \eqref{eq:CLTLn2} in the LCS case, when $\mathbf{X}_{n}$ and $\mathbf{Y}_{n}$ are independent i.i.d. random sequences such that $\text{Var}(L_{n})$ admits a sublinear (in $n$) lower bound. The result in~\cite{HoudreIslak:2015} upper-bounds the Monge-Kantorovich-Wasserstein distance, and in turn the Kolmogorov distance, implying weak convergence. The limiting theorem in~\cite{HoudreIslak:2015} was extended to multiple independent i.i.d. sequences with a general score function in~\cite{GongHoudreIslak:2016}. In contrast to~\cite{HoudreIslak:2015}, the result in~\cite{GongHoudreIslak:2016} directly upper-bounds the Kolmogorov distance, which improves the rate of weak convergence, but still requires a (looser) sublinear variance lower bound. Hence, establishing a linear-order variance lower bound for a general distribution of $X_{1}$ implies immediately the corresponding limiting theorem \eqref{eq:CLTLn2}. Therefore, there is a direct connection between a variance lower bound and the limiting theorem \eqref{eq:CLTLn2}, leading to an extra motivation to study the rate of convergence of $\text{Var}(L_{n})$ as well as all the other central moments of $L_{n}$.

The study of the asymptotic order of $\text{Var}(L_{n})$ in the case of LCS was first proposed by Chv\'{a}tal and Sankoff~\cite{ChvatalSankoff:1975}, who conjectured that $\text{Var}(L_{n})=o(n^{2/3})$, for $\mathbf{X}_{n}$ and $\mathbf{Y}_{n}$ independent i.i.d. symmetric Bernoulli sequences, based on Monte-Carlo simulations. By the Efron-Stein inequality, in the case of independent i.i.d. sequences with score functions satisfying \eqref{eq:MaxScoreChange} below, there exist a universal constant (independent of $n$) $C_{2}>0$, such that
\begin{align}\label{eq:UpperBoundLn}
\text{Var}(L_{n})\leq C_{2}\,n,\quad\text{for all }\,n\in\mathbb{N}.
\end{align}
(see Section \ref{sec:upper} for more details). For the LCS case, this result was proved by Steele~\cite{Steele:1986}. In~\cite{Waterman:1994}, Waterman asked whether or not the linear bound on the variance can be improved, at least in the LCS case. His simulations showed that, in some special cases (including the LCS case), $\text{Var}(L_{n})$ should grow linearly in $n$. These simulations suggest a linear lower bound on $\text{Var}(L_{n})$, which would invalidate the conjecture of Chv\'{a}tal and Sankoff. In the past ten years, the asymptotic behavior of $\text{Var}(L_{n})$ has been investigated under various choices of sequences $\mathbf{X}_{n}$ and $\mathbf{Y}_{n}$ (cf.~\cite{BonettoMatzinger:2006}, \cite{DurringerLemberMatzinger:2007}, \cite{HoudreMatzinger:2007}, \cite{LemberMatzinger:2009}, \cite{LemberMatzingerTorres:2012(2)}, \cite{AmsaluHoudreMatzinger:2014}, \cite{AmsaluHoudreMatzinger:2016}, etc). In particular, in~\cite{LemberMatzingerTorres:2012(2)}, the asymptotic behavior of $\text{Var}(L_{n})$ is investigated within two scenarios. In the first one, a general scoring function with a high gap penalty is considered, and $\mathbf{X}_{n}$ and $\mathbf{Y}_{n}$ are independent i.i.d sequences drawn from a finite alphabet (see also Theorem \ref{thm:app1} and Remark \ref{rem:GapPrice} below). In the second one, the LCS case is studied when both $\mathbf{X}_{n}$ and $\mathbf{Y}_{n}$ are binary sequence having a multinomial block structure, which is a generalization of the model in~\cite{Torres:2009} (the first analysis of the so-called {\it high entropy} case). In many ways, the present paper complements~\cite{LemberMatzingerTorres:2012(2)}. The results of~\cite{LemberMatzinger:2009} (the so-called {\it low entropy} case) were generalized in~\cite{HoudreMa:2016}, and a linear variance lower bound was proved for an arbitrary finite alphabet (not just binary) and for the central moments of arbitrary order (not just the variance), but still under a strongly asymmetric distribution over the letters.

The goal of the present paper is to study in this general framework the lower bounds on generalized moments of the optimal alignment score $L_{n}$. To start with, in the next section we quickly show how to obtain upper moments and exponential bounds on $L_{n}$. In Section \ref{sec:GenMomentBounds}, we develop a general argument to find lower bounds for the generalized central moments of $L_{n}$. These lower bounds are obtained under a main assumption, and sufficient conditions for the validity of this main assumption are provided. Next, we relax some of these conditions (so that the main assumption is no longer satisfied), and develop a ``uniform approximation" method giving rougher lower bounds for the generalized central moments of $L_{n}$. In Section \ref{sec:LCS}, the general results are applied to the binary LCS case, and some refinements of those lower bounds are also presented. Finally, in Section \ref{sec:Rate}, we obtain a lower bound on some moment generating function, which is then used to derive a (local) upper bound on the rate function associated with the binary LCS case.

\section{Upper Bounds on Moments and Exponential Moments}\label{sec:upper}

In this short section, we briefly review the upper bounds on central moments and exponential moments of $L_{n}$. Recall that for any constant gap price $\delta\in\mathbb{R}$, changing the value of one of the $2n$ random variables $X_{1},\ldots,X_{n}$ and $Y_{1},\ldots,Y_{n}$ changes the value of $L_{n}$ by at most $K$, where
\begin{align}\label{eq:MaxScoreChange}
K:=\max_{u,v,w\in\mathcal{A}}|S(u,v)-S(u,w)|.
\end{align}
Thus, by Hoeffding's exponential martingale inequality, for every $t>0$,
\begin{align}\label{eq:HoeffdingIneqOptScore}
\mathbb{P}\left(\left|L_{n}-\mathbb{E}(L_{n})\right|\geq t\right)\leq 2\,e^{-t^{2}/(nK^{2})},
\end{align}
which indicates that $V_{n}:=(L_{n}-\mathbb{E}(L_{n}))/\sqrt{n}$ is subgaussian. Hence, for any $r>0$,
\begin{align*}
\mathbb{E}\left(\left|V_{n}\right|^{r}\right)\leq r\,2^{r/2}\left(9K^{2}\right)^{r/2}\Gamma\left(\frac{r}{2}\right)=r(18)^{r/2}\,\Gamma\left(\frac{r}{2}\right)K^{r}.
\end{align*}

The above inequality provides an upper estimate, of the correct order, on the central absolute moments of $L_{n}$. Some further simple and direct computations allow to improve the constant. Let $r\geq 2$, $x>0$, and let $\widetilde{V}_{n}:=|L_{n}-\mathbb{E}(L_{n})|$. From \eqref{eq:HoeffdingIneqOptScore},
\begin{align*}
\mathbb{E}\left(\widetilde{V}_{n}^{r}\right)=\int_{0}^{\infty}\mathbb{P}\left(\widetilde{V}_{n}\geq t^{1/r}\right)dt\leq x+2\int_{x}^{\infty}e^{-t^{2/r}/(nK^{2})}\,dt.
\end{align*}
Minimizing in $x$, i.e., taking $x=\left(K^{2}(\ln 2)n\right)^{{r/2}}$, and changing variables $u=t^{2/r}/(K^{2}n)$, lead to
\begin{align*}
\mathbb{E}\!\left(\widetilde{V}_{n}^{r}\right)\leq K^{r}\left[(\ln 2)^{r/2}+r\int_{\ln 2}^{\infty}e^{-u}u^{r/2-1}du\right]n^{r/2}=:C(r)\,n^{r/2}.
\end{align*}
When $x=0$, the corresponding constant is slightly bigger than $C(r)$, and is given by
\begin{align*}
D(r):=rK^{r}\int_{0}^{\infty}e^{-u}u^{r/2-1}du=rK^{r}\,\Gamma\left(\frac{r}{2}\right).
\end{align*}
\begin{remark}\label{rem:LCSUpConst}
For the LCS case, the following constant was obtained, as a consequence of a general tensorization inequality, in~\cite{HoudreMa:2016}:
\begin{align*}
E(r):=(r-1)^{r}2^{r/2-1}\mathbb{P}\left(X_{1}\neq Y_{1}\right).
\end{align*}
For $r$ close to 2, $E(r)$ is smaller than $C(r)$. Indeed,
\begin{align*}
E(2)=\mathbb{P}\left(X_{1}\neq Y_{1}\right)\leq 1,\quad C(2)=1+\ln 2,\quad D(2)=2.
\end{align*}
However, for $r>2$ large, $E(r)$ might be bigger, e.g.,
\begin{align*}
E(4)=162\,\mathbb{P}\left(X_{1}\neq Y_{1}\right),\quad D(4)=4\Gamma(2)=4.
\end{align*}
\end{remark}

In addition to moment estimates, the above methodology also provides an upper bound on the moment generating function of $|V_{n}|=\widetilde{V}_{n}/\sqrt{n}=|L_{n}-\mathbb{E}(L_{n})|/\sqrt{n}$. Let $t>0$ and $r>1$. Then, with $K=1$ and using \eqref{eq:HoeffdingIneqOptScore},
\begin{align*}
\mathbb{P}\left(e^{t|V_{n}|}\geq r\right)=\mathbb{P}\left(\widetilde{V}_{n}\geq\frac{\sqrt{n}\ln r}{t}\right)\leq 2\,e^{-\left(\ln r\right)^{2}/t^{2}}.
\end{align*}
Thus,
\begin{align}
M(t):=\mathbb{E}\left(e^{t\left|V_{n}\right|}\right)&\leq 1+2\int_{1}^{\infty}e^{-\left(\ln r\right)^{2}/t^{2}}\,dr=1+2\int_{0}^{\infty}\exp\left(-\frac{x^{2}}{t^{2}}+x\right)dx\nonumber\\
\label{eq:UpperBoundMGFVn} &\leq 1+t\sqrt{\pi}\left(1+\text{erf}\left(\frac{t}{2}\right)\right)e^{t^{2}/4}<1+2t\sqrt{\pi}\,e^{t^{2}/4},
\end{align}
where $\text{erf}\,(\cdot)$ denotes the Gaussian error function.
\begin{remark}\label{rem:UpperBoundSerExp}
The upper bound on the moment generating function can clearly also be found through the upper bounds on the moments. Indeed, when $K=1$ and with the bound
\begin{align*}
\mathbb{E}\left(\left|V_{n}\right|^{r}\right)\leq D(r)=r\,\Gamma\left(\frac{r}{2}\right)=2\,\Gamma\left(\frac{r+2}{2}\right)=2\sqrt{\pi}\,\mathbb{E}\left(\left|\xi\right|^{r+1}\right),\quad r\in\mathbb{N},
\end{align*}
where $\xi$ is a $\mathcal{N}(0,1/2)$ random variable, it follows that
\begin{align*}
\mathbb{E}\left(e^{t|V_{n}|}\right)=1+\sum_{r=1}^{\infty}\mathbb{E}\left(\left|V_{n}\right|^{r}\right)\frac{t^{r}}{r!}\leq 1+2\sqrt{\pi}\sum_{r=1}^{\infty}\mathbb{E}\left(\left|\xi\right|^{r+1}\right)\frac{t^{r}}{r!}=1+\frac{2\sqrt{\pi}}{t}\sum_{k=2}^{\infty}k\,\mathbb{E}\left(\left|\xi\right|^{k}\right)\frac{t^{k}}{k!}.
\end{align*}
Since,
\begin{align*}
\mathbb{E}\left(\left|\xi\right|e^{t|\xi|}\right)=\mathbb{E}\left(|\xi|\right)+\frac{1}{t}\sum_{k=2}^{\infty}k\,\mathbb{E}\left(\left|\xi\right|^{k}\right)\frac{t^{k}}{k!}=\frac{1}{\sqrt{\pi}}+\frac{1}{t}\sum_{k=2}^{\infty}k\,\mathbb{E}\left(\left|\xi\right|^{k}\right)\frac{t^{k}}{k!},
\end{align*}
while,
\begin{align*}
\mathbb{E}\left(\left|\xi\right|e^{t|\xi|}\right)=\frac{d}{dt}\mathbb{E}\left(e^{t|\xi|}\right)=\frac{d}{dt}\left[\left(1+\emph{erf}\left(\frac{t}{2}\right)\right)e^{t^{2}/4}\right]=\frac{t}{2}\left(1+\emph{erf}\left(\frac{t}{2}\right)\right)e^{t^{2}/4}+\frac{1}{\sqrt{\pi}},
\end{align*}
the upper bound in \eqref{eq:UpperBoundMGFVn} is recovered:
\begin{align*}
\mathbb{E}\left(e^{t|V_{n}|}\right)\leq 1+2\sqrt{\pi}\left(\mathbb{E}\left(\left|\xi\right|e^{t|\xi|}\right)-\mathbb{E}\left(|\xi|\right)\right)=1+t\sqrt{\pi}\left(1+\text{erf}\left(\frac{t}{2}\right)\right)e^{t^{2}/4}.
\end{align*}
\end{remark}

\section{Lower Bounds on Generalized Moments}\label{sec:GenMomentBounds}

Let $\Phi:\mathbb{R}^{+}\rightarrow\mathbb{R}^{+}$ be a convex, nondecreasing function. The main objective of this section is to find a lower bound for the expectation
\begin{align}\label{mir}
\mathbb{E}\left(\Phi\left(\left|L_{n}\!\left(\mathbf{X}_{n};\mathbf{Y}_{n}\right)-\mathbb{E}\left(L_{n}\!\left(\mathbf{X}_{n};\mathbf{Y}_{n}\right)\right)\right|\right)\right),
\end{align}
where, again
\begin{align*}
L_{n}=L_{n}\!\left(\mathbf{X}_{n};\mathbf{Y}_{n}\right)=\max_{(\boldsymbol\pi,\boldsymbol\mu)}U_{(\boldsymbol\pi,\boldsymbol\mu)}\!\left(\mathbf{X}_{n};\mathbf{Y}_{n}\right),
\end{align*}
and where again the maximum is taken over all possible alignments.
\begin{remark}
In what follows, we deal with lower bounds on \eqref{mir} with no further restrictions on $\Phi$. This is somehow different from the upper bound case, where to the best of our knowledge, there is no uniform approach that applies to every convex nondecreasing $\Phi$. Indeed, in the previous section, we established upper bounds for $\Phi(x)=|x|^{r}$, $r>0$, and $\Phi(x)=e^{tx}$, $t>0$. As shown in~\cite{HoudreMa:2016}, Efron-Stein type inequalities can be applied to $\Phi(x)=|x|^{r}$, $r>0$, via the Burkholder's inequality, giving the aforementioned constant $E(r)$. Another generalization of the Efron-Stein inequality, given, for example, in~\cite[Chapter 14]{BoucheronLugosiMassart:2013}, is the subadditivity inequality of the so called $\Phi$-Entropy. However, it only applies to those functions $\Phi$ having strictly positive second derivative and such that $1/\Phi''$ is concave, e.g., $\Phi(x)=|x|^{r}$, for $r\in(1,2]$, or $\Phi(x)=x\log x$. Moreover, the upper bounds on \eqref{mir} for those $\Phi$ are linear in $n$ as easily seen via, e.g., \cite[Theorem 14.6]{BoucheronLugosiMassart:2013}, together with \eqref{eq:MaxScoreChange}. In~\cite[Chapter 15]{BoucheronLugosiMassart:2013}, some further generalizations of Efron-stein inequality for $\Phi(x)=|x|^{r}$, $r\geq 2$, are also obtained, which, together with \eqref{eq:MaxScoreChange}, also imply upper bounds of order $n^{r/2}$ on the $r$-th central moments (cf. Theorem 15.5 therein).
\end{remark}

\subsection{The Random Variable $U_{n}$ and the Set $\mathcal{U}_{n}$}

Recall that $\mathcal{X}_{n}$ is the sample space of $\mathbf{X}_{n}$ and $\mathbf{Y}_{n}$, so that $\mathcal{X}_{n}\times\mathcal{X}_{n}$ is the sample space of $\mathbf{Z}_{n}=(\mathbf{X}_{n};\mathbf{Y}_{n})$. In the sequel, we consider a function $\mathfrak{u}:\,\mathcal{X}_{n}\times\mathcal{X}_{n}\rightarrow\mathbb{Z}$, so that $U_{n}:=\mathfrak{u}(\mathbf{Z}_{n})$ is an integer-valued random variable. Let $\mathcal{S}_{n}$ be the support of $U_{n}$, and for every $u\in\mathcal{S}_{n}$, let
\begin{align*}
\ell_{n}(u):=\mathbb{E}\left(L_{n}(\mathbf{Z}_{n})\,|\,U_{n}=u\right),
\end{align*}
and set $\mu_{n}:=\mathbb{E}(L_{n}(\mathbf{Z}_{n}))$. Since $x\mapsto\Phi(|x-\mu_{n}|)$ is convex,
\begin{align*}
\mathbb{E}\left(\Phi\left(|L_{n}(\mathbf{Z}_{n})-\mu_{n}|\right)\,|\,U_{n}\right)\geq\Phi\left(\left|\mathbb{E}\left(L_{n}(\mathbf{Z}_{n})|\,U_{n}\right)-\mu_{n}\right|\right)=\Phi\left(\left|\ell_{n}(U_{n})-\mu_{n}\right|\right),
\end{align*}
so that
\begin{align*}
\mathbb{E}\left(\Phi\left(\left|L_{n}(\mathbf{Z}_{n})-\mu_{n}\right|\right)\right)\geq\mathbb{E}\left(\Phi\left(\left|\ell_{n}(U_{n})-\mu_{n}\right|\right)\right).
\end{align*}

Let $\delta>0$, and let the set $\mathcal{U}_{n}\subset\mathcal{S}_{n}$ be such that
\begin{align}\label{eq:MainAssump}
\ell_{n}(u_{2})-\ell_{n}(u_{1})\geq\delta,\quad\text{for any }\,u_{1},u_{2}\in\mathcal{U}_{n}\,\text{ with }\,u_{1}<u_{2}.
\end{align}
Since $\mathcal{S}_{n}$ is finite, so is $\mathcal{U}_{n}$. Therefore, set
\begin{align}\label{eq:DefUnk0}
\mathcal{U}_{n}:=\{u_{1},\ldots,u_{m}\}\quad\text{and}\quad k_{0}:=\max_{i=1,\,\cdots,m-1}\left(u_{i+1}-u_{i}\right),
\end{align}
where $u_{1}<u_{2}<\cdots<u_{m}$ and where $m=m(n)$ depends on $n$. Clearly, except perhaps in some very special cases, such a set $\mathcal{U}_{n}$ formally always exists (it might even consists of two elements). As we will see, $\mathcal{U}_{n}$ becomes useful if $\mathbb{P}(U_{n}\in\mathcal{U}_{n})\geq c_{0}$, for $\delta$ and $c_{0}$ independent of $n$.
\begin{lemma}\label{lemma:MeanValelln}
Under \eqref{eq:MainAssump}, for every $u_{i},u_{j}\in\mathcal{U}_{n}$,
\begin{align*}
\left|\ell_{n}(u_{i})-\ell_{n}(u_{j})\right|\geq\frac{\delta}{k_{0}}\left|u_{i}-u_{j}\right|.
\end{align*}
Therefore, there exists $a_{n}\in[u_{1},u_{m}]$, such that for every $u_{i}\in\mathcal{U}_{n}$,
\begin{align}\label{eq:MainAssumpA}
\left|\ell_{n}(u_{i})-\mu_{n}\right|\geq\frac{\delta}{k_{0}}\left|u_{i}-a_{n}\right|.
\end{align}
\end{lemma}

\noindent
\textbf{Proof.} Let $g:[u_{1},u_{m}]\rightarrow\mathbb{R}^{+}$ be a differentiable function, such that $g(u_{i})=\ell_{n}(u_{i})$, for every $i=1,\ldots,m$, and such that $g'(x)\geq\delta/k_{0}$, for every $x\in(u_{1},u_{m})$. Then, one of the following three mutually exclusive possibilities holds:
\begin{itemize}
\item $\mu_{n}<\ell_{n}(u_{1})$,
\item $\ell_{n}(u_{1})\leq\mu_{n}\leq\ell_{n}(u_{m})$,
\item $\ell_{n}(u_{m})<\mu_{n}$.
\end{itemize}
In the first case, for every $u_{i}\in\mathcal{U}_{n}$, $i=1,\ldots,m$,
\begin{align*}
\left|\ell_{n}(u_{i})-\mu_{n}\right|=\ell_{n}(u_{i})-\mu_{n}\geq\ell_{n}(u_{i})-\ell_{n}(u_{1})\geq\frac{\delta}{k_{0}}\left(u_{i}-u_{1}\right),
\end{align*}
so \eqref{eq:MainAssumpA} holds with $a_{n}=u_{1}$. The third case can be handled similarly with $a_{n}=u_{m}$. For the middle case, since $g$ is increasing and continuous, there exists $a_{n}\in[u_{1},u_{m}]$ such that $\mu_{n}=g(a_n)$. Hence, by the mean value theorem,
\begin{align*}
\left|\ell_{n}(u_{i})-\mu_{n}\right|=\left|g(u_{i})-g(a_{n})\right|\geq\frac{\delta}{k_{0}}\left|u_{i}-a_{n}\right|,
\end{align*}
completing the proof.\hfill $\Box$

\medskip
Now, in view of \eqref{eq:MainAssumpA}, since $\Phi$ is monotone and non-negative,
\begin{align*}
\mathbb{E}\left(\Phi\left(\left|\ell_{n}(U_{n})-\mu_{n}\right|\right)\right)\geq\sum_{u_{i}\in\mathcal{U}_{n}}\Phi\left(\left|\ell_{n}(u_{i})-\mu_{n}\right|\right)\mathbb{P}\left(U_{n}=u_{i}\right)\geq\sum_{u_{i}\in\mathcal{U}_{n}}\Phi\left(\frac{\delta}{k_{0}}\left|u_{i}-a_{n}\right|\right)\mathbb{P}\left(U_{n}=u_{i}\right),
\end{align*}
i.e., the lower bound on the generalized moment is of the form
\begin{align}\label{eq:BasicEstEPhi}
\mathbb{E}\!\left(\Phi\!\left(\left|L_{n}(\mathbf{Z}_{n})-\mu_{n}\right|\right)\right)\geq\mathbb{E}\!\left(\Phi\!\left(\left|\ell_{n}(U_{n})-\mu_{n}\right|\right)\right)&\geq\sum_{u_{i}\in\mathcal{U}_{n}}\Phi\left(\frac{\delta}{k_{0}}\left|u_{i}-a_{n}\right|\right)\mathbb{P}\left(U_{n}=u_{i}\right)\\
\label{eq:BasicEstEPhi2} &=\mathbb{E}\!\left(\!\left.\Phi\!\left(\frac{\delta}{k_{0}}\left|U_{n}-a_{n}\right|\right)\,\right|U_{n}\in\mathcal{U}_{n}\right)\mathbb{P}\!\left(U_{n}\in\mathcal{U}_{n}\right).
\end{align}
When $\Phi(x)=|x|^{r}$, for some $r\geq 1$, then \eqref{eq:BasicEstEPhi2} becomes
\begin{align*}
\mathbb{E}\left(\left|L_{n}(\mathbf{Z}_{n})-\mu_{n}\right|^{r}\right)\geq\left(\frac{\delta}{k_{0}}\right)^{r}\mathbb{E}\left(\left.\left|U_{n}-a_{n}\right|^{r}\,\right|U_{n}\in\mathcal{U}_{n}\right)\mathbb{P}\left(U_{n}\in\mathcal{U}_{n}\right).
\end{align*}
If $\mathbb{P}(U_{n}\in\mathcal{U}_{n})\geq c_{0}>0$, with $\delta$, $k_{0}$ and $c_{0}$ all independent of $n$, then
\begin{align*}
\mathbb{E}\left(\left|L_{n}(\mathbf{Z}_{n})-\mu_{n}\right|^{r}\right)\geq\left(\frac{\delta}{k_{0}}\right)^{r}c_{0}\,\mathbb{E}\left(\left.\left|U_{n}-a_{n}\right|^{r}\,\right|U_{n}\in\mathcal{U}_{n}\right)\geq\left(\frac{\delta}{k_{0}}\right)^{r}c_{0}\,\min_{a\in\mathbb{R}}\mathbb{E}\left(\left.\left|U_{n}-a\right|^{r}\,\right|U_{n}\in\mathcal{U}_{n}\right).
\end{align*}
Thus a lower bound on the centered absolute $r$th moment of $L_{n}(\mathbf{Z}_{n})$ is obtained via conditioning on $U_{n}$. Typically, in applications, the random variable $U_{n}$ counts the number of some particular letters in $\mathbf{Z}_{n}$, and, therefore, it has a binomial distribution. In this case, the right hand side of \eqref{eq:BasicEstEPhi2} is relatively easy to compute.
\begin{remark}\label{rem:Phin}
Note that the argument leading to \eqref{eq:BasicEstEPhi2} is independent of $n$. Therefore, \eqref{eq:BasicEstEPhi} continues to hold whenever $\Phi$ depends on $n$, provided $\Phi_{n}$ remains convex and non-decreasing.
\end{remark}

\subsection{On the Existence of $U_{n}$ and $\mathcal{U}_{n}$}\label{subsec:ExistUns}

As shown in the previous subsection, the lower bound on the $\Phi$-moment can be obtained via the random variable $U_{n}$, provided that there exists a universal constant $\delta>0$ such that the corresponding set $\mathcal{U}_{n}$ has a large enough probability. In what follows, we show that the existence of a suitable $\mathcal{U}_{n}$ can be guaranteed using a {\it random transformation},
\begin{align*}
\mathcal{R}:\Omega\times\mathcal{X}_{n}\times\mathcal{X}_{n}\rightarrow\mathcal{X}_{n}\times\mathcal{X}_{n},
\end{align*}
such that, for most of the outcomes of $Z$, $\mathcal{R}$ increases the score by at least some fixed amount $\varepsilon_{0}>0$. Here, with a slight abuse of notation, $\Omega$ is a sample space for the randomness involved in $\mathcal{R}$. More precisely, the random transformation should be
such that there exists a set $B_{n}\subset\mathcal{X}_{n}\times\mathcal{X}_{n}$ having a relatively large probability so that, for every $\mathbf{z}_{n}\in B_{n}$, the expected score of $\mathcal{R}(\mathbf{z}_{n})$ exceeds the score of $\mathbf{z}_{n}$ by $\varepsilon_{0}$, namely,
\begin{align*}
\mathbb{E}\left(L_{n}(\mathcal{R}(\mathbf{z}_{n}))\right)\geq L_{n}(\mathbf{z}_{n})+\varepsilon_{0}.
\end{align*}
Above $\mathbb{E}$ denotes the expectation over the randomness involved in $\mathcal{R}$ (i.e., the integral over $\Omega$). The requirement that the set $B_{n}$ has a relatively large probability is formalized by requiring that
\begin{align}\label{eq:ReqRanTrans}
\Delta_{n}(\varepsilon_0):=\mathbb{P}\left(\mathbb{E}\left(\left.L_{n}(\mathcal{R}(\mathbf{Z}_{n}))-L_{n}(\mathbf{Z}_{n})\,\right|\mathbf{Z}_{n}\right)<\varepsilon_{0}\right)\rightarrow 0,\quad\text{as }\,n\rightarrow\infty.
\end{align}
We will see that the faster $\Delta_{n}(\varepsilon_{0})$ tends to zero, the better the lower bound is.

The random transformation should be associated with $U_{n}$ and the set $\mathcal{U}_{n}=\{u_{1},\ldots,u_{m}\}$ so that when $U_{n}$ takes its values in $\mathcal{U}_{n}$, then $\mathcal{R}$ increases its value in $\mathcal{U}_{n}$ and has no effect on the conditional distribution of $\mathbf{Z}_{n}$. To formalize this requirement, for every $u\in\mathcal{U}_{n}$, let $\mathbb{P}^{(u)}$ be the law of $\mathbf{Z}_{n}$ given $U_{n}=u$, i.e.,
\begin{align*}
\mathbb{P}^{(u)}(A)=\mathbb{P}\left(\left.\mathbf{Z}_{n}\in A\,\right|U_{n}=u\right),\quad A\subset\mathcal{X}_{n}\times\mathcal{X}_{n}.
\end{align*}
Now, for every $u_{i},u_{i+1}\in\mathcal{U}_{n}$, the following implication should hold: if $\mathbf{Z}_{n}^{(u_{i})}$ is a random vector with law $\mathbb{P}^{(u_{i})}$, then
\begin{align}\label{eq:RanTransImpli}
\mathcal{R}\left(\mathbf{Z}_{n}^{(u_{i})}\right)\sim\mathbb{P}^{(u_{i+1})},
\end{align}
which means that for every $\mathbf{z}_{n}\in\mathcal{X}_{n}\times\mathcal{X}_{n}$,
\begin{align*}
\mathbb{P}\left(\left.\mathcal{R}(\mathbf{Z}_{n})=\mathbf{z}_{n}\,\right|U_{n}=u_{i}\right)=\mathbb{P}\left(\left.\mathbf{Z}_{n}=\mathbf{z}_{n}\,\right|U_{n}=u_{i+1}\right).
\end{align*}
In particular, \eqref{eq:RanTransImpli} implies that if $\mathfrak{u}(\mathbf{Z}_{n})=u_{i}$, then $\mathfrak{u}(\mathcal{R}(\mathbf{Z}_{n}))=u_{i+1}$, but the converse statement might not be true.

\medskip
\noindent
\textbf{Fast Convergence Condition.} Our first general lower bound theorem assumes the existence of $\varepsilon_{0}$ so that $\Delta_{n}(\varepsilon_{0})$ converges to zero sufficiently fast. To simplify the notation, in what follows, let $\widetilde{\mathbf{Z}}_{n}:=\mathcal{R}(\mathbf{Z}_{n})$. The following theorem is a generalization of~\cite[Theorem 2.2]{LemberMatzingerTorres:2012(2)}.
\begin{theorem}\label{thm:FastConvThm}
For $n\in\mathbb{N}$, let there exist a random variable $U_{n}=\mathfrak{u}(\mathbf{Z}_{n})$, a set $\mathcal{U}_{n}=\{u_{1},\ldots,u_{m(n)}\}$, a random transformation $\mathcal{R}$ satisfying \eqref{eq:RanTransImpli}, such that the following two assumptions hold:
\begin{itemize}
\item [(i)]  There exists a universal constant $A>0$ such that
    \begin{align}\label{eq:DiffLtildeZLZ}
    L_{n}\big(\widetilde{\mathbf{Z}}_{n}\big)-L_{n}(\mathbf{Z}_{n})\geq -A.
    \end{align}
\item [(ii)] There exists a universal constant $\varepsilon_{0}>0$ such that
    \begin{align}\label{eq:DeltafastVarphi}
    \Delta_{n}(\varepsilon_{0})=o(\varphi(n)),\quad n\rightarrow\infty,
    \end{align}
    where $\Delta_{n}(\varepsilon_{0})$ is given by \eqref{eq:ReqRanTrans}, and where $\varphi$ be a function on $\mathbb{N}$ such that
    \begin{align}\label{eq:GenLowerBoundProbU}
    \mathbb{P}\left(U_{n}=u\right)\geq\varphi(n),\quad\text{for all }\,u\in\mathcal{U}_{n}.
    \end{align}
\end{itemize}
Then, there exists $a_{n}\in[u_{1},u_{m}]$ so that for $n$ large enough,
\begin{align}\label{eq:trm31}
\mathbb{E}\left(\Phi\left(\left|L_{n}(\mathbf{Z}_{n})-\mu_{n}\right|\right)\right)\geq\sum_{u_{i}\in\mathcal{U}_{n}}\Phi\left(\frac{\delta}{k_{0}}\left|u_{i}-a_{n}\right|\right)\mathbb{P}\left(U_{n}=u_{i}\right)\geq\sum_{u_{i}\in\mathcal{U}_{n}}\Phi\left(\frac{\delta}{k_{0}}\left|u_{i}-a_{n}\right|\right)\varphi(n),
\end{align}
where $\delta=\varepsilon_{0}/2$ and $k_{0}(n)=\max_{i=1,\cdots,m-1}(u_{i+1}-u_{i})$.
\end{theorem}
\noindent
\textbf{Proof.} Let $u_{i}\in\mathcal{U}_{n}$, and let $\mathbf{Z}_{n}^{(u_{i})}$ be a random vector with law $\mathbb{P}^{(u_{i})}$. By \eqref{eq:RanTransImpli},
\begin{align*}
\ell_{n}(u_{i+1})=\mathbb{E}\left(L_{n}\!\left(\widetilde{\mathbf{Z}}_{n}^{(u_{i})}\right)\right).
\end{align*}
Hence,
\begin{align*}
\ell_{n}(u_{i+1})-\ell_{n}(u_{i})=\mathbb{E}\left(L_{n}\!\left(\widetilde{\mathbf{Z}}_{n}^{(u_{i})}\right)\right)-\mathbb{E}\left(L_{n}\!\left(\mathbf{Z}_{n}^{(u_{i})}\right)\right)=\mathbb{E}\left[\mathbb{E}\left(\left.L_{n}\!\left(\widetilde{\mathbf{Z}}_{n}^{(u_{i})}\right)-L_{n}\!\left(\mathbf{Z}_{n}^{(u_{i})}\right)\,\right|\mathbf{Z}_{n}^{(u_{i})}\right)\right].
\end{align*}
Let $\varepsilon_{0}$ be as in (ii), and let $B_{n}(\varepsilon_{0})\subset\mathcal{X}_{n}\times\mathcal{X}_{n}$ be the set of outcomes of $\mathbf{Z}_{n}$ such that
\begin{align}\label{eq:DefZBn}
\left\{\mathbb{E}\left(\left.L_{n}\big(\widetilde{\mathbf{Z}}_{n}\big)-L_{n}(\mathbf{Z}_{n})\,\right|\mathbf{Z}_{n}\right)\geq\varepsilon_{0}\right\}=\left\{\mathbf{Z}_{n}\in B_{n}(\varepsilon_{0})\right\}.
\end{align}
By \eqref{eq:DiffLtildeZLZ}, for any pair of sequences $z_{n}$, when applying the transformation $\mathcal{R}$, the score can decrease by at most $-A$. Hence
\begin{align*}
\mathbb{E}\left[\mathbb{E}\left(\left.L_{n}\!\left(\widetilde{\mathbf{Z}}_{n}^{(u_{i})}\right)-L_{n}\!\left(\mathbf{Z}_{n}^{(u_{i})}\right)\,\right|\mathbf{Z}_{n}^{(u_{i})}\right)\right]\geq\varepsilon_{0}\,\mathbb{P}\left(\mathbf{Z}_{n}^{(u_{i})}\in B_{n}(\varepsilon_{0})\right)-A\,\mathbb{P}\left(\mathbf{Z}_{n}^{(u_{i})}\not\in B_{n}(\varepsilon_{0})\right).
\end{align*}
By definition, $\mathbb{P}(\mathbf{Z}_{n}\not\in B_{n}(\varepsilon_{0}))=\Delta_{n}(\varepsilon_{0})$. Therefore, by \eqref{eq:DeltafastVarphi} and \eqref{eq:GenLowerBoundProbU},
\begin{align*}
\mathbb{P}\left(\mathbf{Z}_{n}^{(u_{i})}\not\in B_{n}(\varepsilon_{0})\right)=\mathbb{P}\left(\left.\mathbf{Z}_{n}\not\in B_{n}(\varepsilon_{0})\,\right|U_{n}=u_{i}\right)\leq\frac{\Delta_{n}(\varepsilon_{0})}{\mathbb{P}\left(U_{n}=u_{i}\right)}\leq\frac{\Delta_{n}(\varepsilon_{0})}{\varphi(n)}=o(1),\quad
n\rightarrow\infty.
\end{align*}
Choose $n_{0}$ large enough (depending on $\varepsilon_{0}$ and $A$) so that for any $n>n_{0}$,
\begin{align*}
\varepsilon_{0}\left(1-\frac{\Delta_{n}(\varepsilon_{0})}{\varphi(n)}\right)-A\,\frac{\Delta_{n}(\varepsilon_{0})}{\varphi(n)}\geq\frac{\varepsilon_{0}}{2}.
\end{align*}
Therefore, for any $n>n_{0}$ and any $u_{i}\in\mathcal{U}_{n}$,
\begin{align*}
\ell_{n}(u_{i+1})-\ell_{n}(u_{i})\geq\frac{\varepsilon_{0}}{2}=:\delta.
\end{align*}
This proves \eqref{eq:MainAssump}, which entails \eqref{eq:MainAssumpA} and thus \eqref{eq:trm31}. The proof is now complete.\hfill $\Box$
\begin{remark}\label{rem:SufCondMainAssump}
The assumption \eqref{eq:DiffLtildeZLZ} typically holds for any meaningful function $\mathfrak{u}$. The difficulties in applying Theorem \ref{thm:FastConvThm} lie in finding $\mathfrak{u}$, $\mathcal{R}$ and $\mathcal{U}_{n}$, such that \eqref{eq:DeltafastVarphi} holds for a given $\varepsilon_{0}$. Since typically for every $u\in\mathcal{S}_{n}$, $\mathbb{P}(U_{n}=u)\rightarrow 0$, $\Delta_{n}(\varepsilon_{0})$ has to converge to zero sufficiently fast as $n\rightarrow\infty$. The larger the probability $\mathbb{P}(U_{n}\in\mathcal{U}_n)$ is (recall that it has to be bounded away from zero), the smaller $\varphi(n)$ is, and the faster the convergence to zero of $\Delta_{n}(\varepsilon_{0})$ is.
\end{remark}

\subsection{Uniform Approximation}\label{subsec:UnifApprox}

As mentioned in Remark \ref{rem:SufCondMainAssump}, the most important assumption in Theorem \ref{thm:FastConvThm} is the existence of the random transformation $\mathcal{R}$ and of the corresponding random variable $U_{n}$ such that (for some $\varepsilon_{0}>0$) the probability $\Delta_{n}(\varepsilon_{0})$ converges to zero sufficiently fast. For all purposes, this is often a very specific requirement and so far, the existence of suitable $\mathcal{R}$ and $U_{n}$ has only been shown for very specific models (cf.~\cite{HoudreMatzinger:2007}, \cite{LemberMatzinger:2009}, \cite{LemberMatzingerTorres:2012(2)} and~\cite{HoudreMa:2016}). It is therefore important to relax the assumption of fast convergence, which is the goal of the present subsection. Below, the condition \eqref{eq:DeltafastVarphi} is replaced by an arbitrarily slowly convergent sequence $\Delta_{n}(\varepsilon_{0})\rightarrow 0$.
\begin{theorem}\label{thm:UnifConv}
For every $n\in\mathbb{N}$, let there exist a random variable $U_{n}=\mathfrak{u}(\mathbf{Z}_{n})$, a set $\mathcal{U}_{n}=\{u_{1},\ldots,u_{m(n)}\}\subset\mathcal{S}_{n}$, a random transformation $\mathcal{R}$, a positive function $\varphi$ on $\mathbb{N}$, and a constant $A>0$, such that \eqref{eq:RanTransImpli}, \eqref{eq:DiffLtildeZLZ} and \eqref{eq:GenLowerBoundProbU} are satisfied. Moreover, let there exist a universal constant $c>0$, such that
\begin{align}\label{eq:LowerBoundSizeUn}
m(n)\geq c\,\varphi^{-1}(n),\quad\text{for all }\,n\in\mathbb{N}.
\end{align}
If for an $\varepsilon_{0}>0$, $\Delta_{n}(\varepsilon_{0})\rightarrow 0$, then for $n$ large enough,
\begin{align*}
\mathbb{E}\left(\Phi\left(\left|L_{n}(\mathbf{Z}_{n})-\mu_{n}\right|\right)\right)\geq\sum_{j=1}^{\frac{c}{8}\varphi^{-1}(n)}\Phi\left(\frac{\varepsilon_{0}}{2}\left(\frac{c}{8\varphi(n)}+j\right)\right)\varphi(n)+\frac{c}{8}\,\Phi\left(\frac{\varepsilon_{0}c}{16\varphi(n)}\right)\geq\frac{c}{4}\,\Phi\left(\frac{\varepsilon_{0}c}{16\varphi(n)}\right). \end{align*}
\end{theorem}
\noindent
\textbf{Proof.} Let $\varepsilon_{0}$ be as in the assumption of the theorem and let the set $\mathcal{U}_{n}^{0}\subseteq\mathcal{S}_{n}$ be defined as follows:
\begin{align}\label{eq:DefUn0}
u\in\mathcal{U}_{n}^{0}\quad\Leftrightarrow\quad\mathbb{P}\left(\left.\mathbf{Z}_{n}\not\in B_{n}(\varepsilon_{0})\right|U_{n}=u\right)\leq\sqrt{\Delta_{n}(\varepsilon_{0})},
\end{align}
where $B_{n}(\varepsilon_{0})$ is given by \eqref{eq:DefZBn}. Thus, by the very definition of $\Delta_{n}(\varepsilon_{0})$,
\begin{align*}
\Delta_{n}(\varepsilon_{0})=\mathbb{P}\left(\mathbf{Z}_{n}\not\in B_{n}(\varepsilon_{0})\right)\geq\sum_{u\not\in\mathcal{U}_{n}^{0}}\mathbb{P}\left(\left.\mathbf{Z}_{n}\not\in B_{n}(\varepsilon_{0})\,\right|U_{n}=u\right)\mathbb{P}\left(U_{n}=u\right)\geq\sqrt{\Delta_{n}(\varepsilon_{0})}\,\mathbb{P}\left(U_{n}\not\in\mathcal{U}_{n}^{0}\right),
\end{align*}
which implies that
\begin{align}\label{eq:ProbUnotUn0}
\mathbb{P}\left(U_{n}\not\in\mathcal{U}_{n}^{0}\right)\leq\sqrt{\Delta_{n}(\varepsilon_{0})}.
\end{align}
When $u_{i}\in\mathcal{U}_{n}\cap\mathcal{U}_{n}^{0}$, using an argument as in the proof of Theorem \ref{thm:FastConvThm}, and by \eqref{eq:DiffLtildeZLZ} and \eqref{eq:DefUn0},
\begin{align}\label{eq:DiffellUnUn0}
\ell_{n}(u_{i+1})-\ell_{n}(u_{i})\geq\varepsilon_{0}\left(1-\sqrt{\Delta_{n}(\varepsilon_{0})}\right)-A\sqrt{\Delta_{n}(\varepsilon_{0})}\geq\frac{\varepsilon_{0}}{2},
\end{align}
provided that $n>n_{1}$, for some positive integer $n_{1}$ large enough. When $u_{i}\not\in\mathcal{U}_{n}\cap\mathcal{U}_{n}^{0}$, by \eqref{eq:DiffLtildeZLZ},
\begin{align}\label{eq:DiffEllNotUnUn0}
\ell_{n}(u_{i+1})-\ell_{n}(u_{i})\geq -A.
\end{align}

Assuming $n>n_{1}$, we are interested in the set $\mathcal{U}_{n}\cap\mathcal{U}_{n}^{0}$ which can be represented as the union of disjoint subintervals of $\mathcal{U}_{n}=\{u_{1},\,\cdots,u_{m}\}$, i.e., $\mathcal{U}_{n}\cap\mathcal{U}_{n}^{0}=\cup_{j=1}^{N}I_{j}$, where for $j=1,\ldots,N$,
\begin{align*}
I_{j}=\left\{u_{1}^{(j)},\ldots,u_{r_{j}}^{(j)}\right\},\quad\text{with }\,u_{1}^{(1)}<\cdots<u_{r_{1}}^{(1)}<u_{1}^{(2)}<\cdots<u_{r_{2}}^{(2)}<\cdots<u_{1}^{(N)}<\cdots<u_{r_{N}}^{(N)},
\end{align*}
are disjoint subintervals of $\mathcal{U}_{n}$. Clearly, the number $N$ of such intervals as well as the intervals $I_{1},\ldots,{I_{N}}$ themselves depend on $n$. On each interval $I_{j}$, \eqref{eq:DiffellUnUn0} implies that the function $\ell_{n}$ increases with a slope of at least $\varepsilon_{0}/2$, and moreover, note that if $u_{i}=u_{r_{j}}^{(j)}$, then $\ell_{n}(u_{i+1})\not\in I_{j}$. Note also that the approach in proving Theorem \ref{thm:FastConvThm} largely applied because $\mathcal{U}_{n}$ was a single interval where $\ell_{n}$ increases. In the present, case the set $\mathcal{U}_{n}\cap\mathcal{U}^{0}_{n}$ consists of several intervals and between them the function $\ell_{n}$ does not necessarily increase. At first, it might appear that the approach of the previous subsection would work, when replacing
\begin{align*}
\sum_{u_{i}\in\mathcal{U}_{n}}\Phi\left(\delta\left|u_{i}-a_{n}\right|\right)\mathbb{P}\left(U_{n}=u_{i}\right)
\end{align*}
with
\begin{align*}
\sum_{u_{i}\in\mathcal{U}_{n}\cap\,\mathcal{U}^{0}_{n}}\Phi\left(\delta\left|u_{i}-a_{n}\right|\right)\mathbb{P}\left(U_{n}=u_{i}\right)=\sum_{j=1}^{N}\sum_{u_{i}\in I_{j}}\Phi\left(\delta\left|u_{i}-a_{n}\right|\right)\mathbb{P}\left(U_{n}=u_{i}\right).
\end{align*}
Unfortunately, $a_{n}$ would now depend on $I_{j}$, and we would have to deal with the sums
\begin{align*}
\sum_{j=1}^{N}\sum_{u_{i}\in I_{j}}\Phi\left(\delta\left|u_{i}-a_{n}^{(j)}\right|\right)\mathbb{P}\left(U_{n}=u_{i}\right).
\end{align*}

To bypass this problem, we proceed differently. Consider the sets
\begin{align*}
J_{j}:=\left\{\ell_{n}\!\left(u_{1}^{(j)}\right),\,\ldots\,,\ell_{n}\!\left(u_{r_{j}}^{(j)}\right)\right\},\quad j=1,\ldots,N,
\end{align*}
i.e., $J_{j}$ is the image of $I_{j}$ under $\ell_{n}$. By \eqref{eq:DiffellUnUn0}, all elements of $J_{j}$ are at least $\varepsilon_{0}/2$-apart from each other (the intervals $J_{j}$ might overlap, although the intervals $I_{j}$ do not). By \eqref{eq:DiffEllNotUnUn0},
\begin{align*}
\sum_{u_{i}\in\mathcal{U}_{n}\setminus\mathcal{U}_{n}^{0}}\left(\ell_{n}(u_{i+1})-\ell_{n}(u_{i})\right)\geq -A\cdot\text{Card}(\mathcal{U}_{n}\setminus\mathcal{U}_{n}^{0})=: -A\,m_{0}(n),
\end{align*}
implying that the sum of the lengths of the (integer) intervals $J_{j}$ differs from the length of $J:=\cup_{j=1}^{N}J_{j}$ by at most $A\,m_{0}(n)$. Formally,
\begin{align*}
\sum_{j=1}^{N}\left(\ell_{n}\!\left(u_{r_{j}}^{(j)}\right)-\ell_{n}\!\left(u_{1}^{(j)}\right)\right)-|J|\leq A\,m_{0}(n),
\end{align*}
where $|J|$ denotes the length of $J$, i.e., the difference between the largest element and the smallest element of $J$, and
$\ell_{n}(u_{r_{j}}^{(j)})-\ell_{n}(u_{1}^{(j)})$ is the length of each $J_{j}$. The number of elements $\varepsilon_{0}/2$-apart needed for covering a (real) interval with length $A\,m_{0}(n)$ is at most $2A\,m_{0}(n)/\varepsilon_{0}+1$. Hence, at most $2A\,m_{0}(n)/\varepsilon_{0}+1$ elements $\varepsilon_{0}/2$-apart will be lost due to overlappings, which implies that the set $J$ contains at least
\begin{align*}
\left|\mathcal{U}_{n}\right|-\frac{2A\,m_{0}(n)}{\varepsilon_{0}}-1=m(n)-\frac{2A\,m_{0}(n)}{\varepsilon_{0}}-1
\end{align*}
elements which are (at least) $\varepsilon_{0}/2$-apart from each other. Using \eqref{eq:GenLowerBoundProbU} and \eqref{eq:ProbUnotUn0},
\begin{align*}
\sqrt{\Delta_{n}(\varepsilon_{0})}\geq\mathbb{P}\left(U_{n}\in\mathcal{U}_{n}\setminus\mathcal{U}^{0}_{n}\right)=\sum_{u\in\mathcal{U}_{n}\setminus\mathcal{U}^{0}_{n}}\mathbb{P}\left(U_{n}=u\right)\geq m_{0}(n)\,\varphi(n),
\end{align*}
and therefore,
\begin{align*}
m_{0}(n)\leq\sqrt{\Delta_{n}(\varepsilon_{0})}\,\varphi^{-1}(n).
\end{align*}
Recalling the assumption \eqref{eq:LowerBoundSizeUn}, there exists $n_{2}\in\mathbb{N}$ such that for all $n>n_{2}$,
\begin{align}\label{eq:LowerBoundNumberEps2Points}
m(n)-\frac{2A\,m_{0}(n)}{\varepsilon_{0}}-1\geq m(n)-\frac{2A\sqrt{\Delta_{n}(\varepsilon_{0})}}{\varepsilon_{0}\varphi(n)}\geq\frac{c\varepsilon_{0}-2A\sqrt{\Delta_{n}(\varepsilon_{0})}}{\varepsilon_{0}\varphi(n)}\geq\frac{c}{2}\varphi^{-1}(n).
\end{align}
Now let
\begin{align*}
\mathcal{B}_{n}:=\left\{u\in\mathcal{U}_{n}:\,\left|\ell_{n}(u)-\mu_{n}\right|\geq\frac{\varepsilon_{0}c}{16\varphi(n)}\right\}=\mathcal{B}_{n}^{+}\cup\mathcal{B}_{n}^{-},
\end{align*}
where
\begin{align*}
\mathcal{B}_{n}^{+}:=\left\{u\in\mathcal{U}_{n}:\,\ell_{n}(u)-\mu_{n}\geq\frac{\varepsilon_{0}c}{16\varphi(n)}\right\}\quad\text{and}\quad\mathcal{B}_{n}^{-}:=\left\{u\in\mathcal{U}_{n}:\,\ell_{n}(u)-\mu_{n}\leq
-\frac{\varepsilon_{0}c}{16\varphi(n)}\right\}.
\end{align*}
Since the interval
\begin{align*}
\left(\mu_{n}-\frac{\varepsilon_{0}c}{16\varphi(n)},\,\mu_{n}+\frac{\varepsilon_{0}c}{16\varphi(n)}\right)
\end{align*}
contains at most $c\varphi^{-1}(n)/4$ elements which are $\varepsilon_{0}/2$-apart from each other, and since the set $J$ contains at least $c\varphi^{-1}(n)/2$ elements which are $\varepsilon_{0}/2$-apart from each other, it follows that $\mathcal{B}_{n}$ contains at least $c\varphi^{-1}(n)/4$ elements (the values of $\ell_{n}$ being $\varepsilon_{0}/2$-apart from each other). This implies that at least one of the two sets $\mathcal{B}_{n}^{+}$ and $\mathcal{B}_{n}^{-}$, say $\mathcal{B}_{n}^{+}$, contains at least $c\varphi^{-1}(n)/8$ elements whose $\ell_{n}$-images are at least $\varepsilon_{0}/2$-apart from each other. Let $\widetilde{\mathcal{B}}_{n}$ be the union of those $c\varphi^{-1}(n)/8$ elements in $\mathcal{B}_{n}^{+}$. Then the set $\mathcal{B}_{n}\setminus\widetilde{\mathcal{B}}_{n}$ consists of at least $c\varphi^{-1}(n)/8$ elements (whose $\ell_{n}$-images are at least $\varepsilon_{0}/2$-apart from each other). Hence, by the monotonicity of $\Phi$,
\begin{align*}
\sum_{i=1}^{m(n)}\Phi\left(\left|\ell_{n}(u_{i})-\mu_{n}\right|\right)\mathbb{P}\left(U_{n}=u_{i}\right)&\geq\sum_{i=1}^{m(n)}\Phi\left(\left|\ell_{n}(u_{i})-\mu_{n}\right|\right)\varphi(n)\geq\sum_{u\in\mathcal{B}_{n}}\Phi\left(\left|\ell_{n}(u)-\mu_{n}\right|\right)\varphi(n)\\
&\geq\sum_{u\in\widetilde{\mathcal{B}}_{n}}\Phi\left(\left|\ell_{n}(u)-\mu_{n}\right|\right)\varphi(n)+\sum_{u\in\mathcal{B}_{n}\setminus\widetilde{\mathcal{B}}_{n}}\Phi\left(\left|\ell_{n}(u)-\mu_{n}\right|\right)\varphi(n)\\
&\geq\sum_{j=1}^{\frac{c}{8}\varphi^{-1}(n)}\!\Phi\!\left(\frac{\varepsilon_{0}}{2}\left(\frac{c}{8\varphi(n)}\!+\!j\right)\right)\!\varphi(n)+\Phi\!\left(\frac{\varepsilon_{0}c}{16\varphi(n)}\right)\frac{c}{8}\varphi^{-1}(n)\varphi(n)\\
&\geq\sum_{j=1}^{\frac{c}{8}\varphi^{-1}(n)}\Phi\left(\frac{\varepsilon_{0}}{2}\left(\frac{c}{8\varphi(n)}+j\right)\right)\varphi(n)+\Phi\left(\frac{\varepsilon_{0}c}{16\varphi(n)}\right)\frac{c}{8}\\
&\geq\frac{c}{4}\,\Phi\left(\frac{\varepsilon_{0}c}{16\varphi(n)}\right).
\end{align*}
The inequality \eqref{eq:BasicEstEPhi} now finishes the proof.\hfill $\Box$
\begin{remark}\label{rem:UniformLargeN}
In both Theorem \ref{thm:FastConvThm} and Theorem \ref{thm:UnifConv}, the thresholds, on large $n$, are independent of the choice of the convex function $\Phi$. In Theorem \ref{thm:FastConvThm}, the threshold $n_{0}$ depends on the sequence $\Delta_{n}(\varepsilon_{0})/\varphi(n)$, but not on $\Phi$. In Theorem \ref{thm:UnifConv}, the threshold on $n$ is obtained in two steps: the first threshold $n_{1}$ is to ensure the validity of \eqref{eq:DiffellUnUn0}, while the second threshold $n_{2}$ is to guarantee the validity of the last inequality in \eqref{eq:LowerBoundNumberEps2Points}. Both steps are independent of the choice of $\Phi$.
\end{remark}
\begin{remark}\label{rem:UnifApprox}
The assumption $m(n)=|\mathcal{U}_{n}|\geq c\varphi^{-1}(n)$ (for all $n\in\mathbb{N}$), together with \eqref{eq:GenLowerBoundProbU}, implies that for every $u\in\mathcal{U}_{n}$, $\mathbb{P}(U_{n}=u)\geq\varphi(n)\geq c\,m^{-1}(n)$. Hence, $\mathcal{U}_{n}$ cannot be too large, since all the atoms in $\mathcal{U}_{n}$ have roughly equal masses, justifying our naming this approach ``uniform approximation". On the other hand, these two assumptions imply that $\mathbb{P}(U\in\mathcal{U}_{n})\geq c$, so that $\mathcal{U}_{n}$ cannot be too small either. In the next section, we will see that the set $\mathcal{U}_{n}$ satisfying the requirements of Theorem \ref{thm:UnifConv} does exist, but its choice must be made more carefully when compared to the one in Theorem \ref{thm:FastConvThm}. This is the price to pay for the arbitrarily slow convergence of $\Delta_{n}(\varepsilon_{0})$.
\end{remark}

\section{The Binary LCS Case}\label{sec:LCS}

As an applications of the general results of the previous section, in the rest of the paper, we consider the binary alphabet $\mathcal{A}=\{0,1\}$, i.e., $\mathbf{X}_{n}$ and $\mathbf{Y}_{n}$ are two independent i.i.d. sequences such that $\mathbb{P}(X_{1}=1)=\mathbb{P}(Y_{1}=1)=p$, $p\in(0,1)$, and $L_{n}(\mathbf{X}_{n};\mathbf{Y}_{n})$ is the length of the LCS of $\mathbf{X}_{n}$ and $\mathbf{Y}_{n}$. Let $\mathfrak{u}(\mathbf{z}_{n})$ count the zeros in $\mathbf{z}_{n}=(\mathbf{x}_{n};\mathbf{y}_{n})$. Thus $U_{n}$ is equal to the number of zeros in $\mathbf{Z}_{n}=(\mathbf{X}_{n};\mathbf{Y}_{n})$ and therefore $U_{n}\sim B(2n,q)$, where $q=1-p$, and its support is $\mathcal{S}_{n}=\{0,1,\ldots,2n\}$.

Let the random transformation $\mathcal{R}$ be defined as follows: given $\mathbf{z}_{n}$, choose randomly (with uniform probability) a one and turn it into a zero. With this transformation, the condition \eqref{eq:RanTransImpli} holds for any $u=0,\ldots,2n-1$. To see this, for any $u=0,\ldots,2n$, let $\mathcal{A}_{n}(u)\subseteq\{0,1\}^{2n}$ consist of all the binary sequences containing exactly $u$ zeros. For any $\mathbf{z}_{n}\in\mathcal{A}_{n}(u)$,
\begin{align*}
\mathbb{P}\left(\left.\mathbf{Z}_{n}=\mathbf{z}_{n}\,\right|U_{n}=u\right)=\frac{\mathbb{P}\left(\mathbf{Z}_{n}=\mathbf{z}_{n},\,U_{n}=u\right)}{\mathbb{P}\left(U_{n}=u\right)}=\frac{\mathbb{P}\left(\mathbf{Z}_{n}=\mathbf{z}_{n}\right)}{\mathbb{P}\left(U_{n}=u\right)}=\frac{q^{u}(1-q)^{2n-u}}{\displaystyle{\binom{u}{2n}}q^{u}(1-q)^{2n-u}}=\binom{u}{2n}^{-1}.
\end{align*}
In other words, $\mathbb{P}_{(u)}$ is the uniform distribution on $\mathcal{A}_{n}(u)$. Now, for any such $u$, let $\mathbf{Z}_{n}^{(u)}$ be a random vector with law $\mathbb{P}^{(u)}$, then $\widetilde{\mathbf{Z}}_{n}^{(u)}=\mathcal{R}(\mathbf{Z}_{n}^{(u)})$ is supported on $\mathcal{A}_{n}(u+1)$. For any $\mathbf{z}_{n}\in\mathcal{A}_{n}(u+1)$, let $0\leq i_{1}<\cdots<i_{u+1}\leq 2n$ be the positions of the zeros in $\mathbf{z}_{n}$, and let $\widehat{\mathbf{z}}_{n}^{(i_{j})}$, $j=1,\ldots,u+1$, be the sequence in $\mathcal{A}_{n}(u)$ obtained by changing the zero of $\mathbf{z}_{n}$ at the position $i_{j}$ to a one. Then,
\begin{align*}
\mathbb{P}\!\left(\widetilde{\mathbf{Z}}_{n}^{(u)}\!=\!\mathbf{z}_{n}\!\right)\!=\!\sum_{j=1}^{u+1}\mathbb{P}\!\left(\left.\widetilde{\mathbf{Z}}_{n}^{(u)}\!=\!\mathbf{z}_{n}\,\right|\mathbf{Z}_{n}^{(u)}\!=\!\widehat{\mathbf{z}}_{n}^{(i_{j})}\right)\mathbb{P}\!\left(\mathbf{Z}_{n}^{(u)}\!=\!\widehat{\mathbf{z}}_{n}^{(i_{j})}\right)\!=\!\frac{u\!+\!1}{2n\!-\!u}\binom{u}{2n}^{-1}\!\!\!\!=\!\binom{u\!+\!1}{2n}^{-1}\!\!\!\!=\!\mathbb{P}^{(u+1)}(\mathbf{z}_{n}).
\end{align*}
Finally note that since we are considering the LCS, any single-entry change in $\mathbf{z}_{n}$ changes the score by at most one. Hence, $A=1$ in \eqref{eq:DiffLtildeZLZ} and
\begin{align}\label{eq:star}
-1\leq\ell_{n}(u+1)-\ell_{n}(u)\leq 1,\quad\text{for any }\,u\in \{0,1,\ldots,2n-1\}.
\end{align}

\subsection{Lower Bounds Under Fast Convergence}\label{subsec:FastConvLCS}

The following theorem is proved in~\cite[Theorem 2.2]{LemberMatzinger:2009} (see also~\cite[Theorem 2.1]{HoudreMa:2016} for a quantitative version with an arbitrary finite alphabet size).
\begin{theorem}\label{thm:LCSLargeProb}
There exist positive constants $\varepsilon_{1}$ and $\varepsilon_{2}$ with $\varepsilon_{1}>\varepsilon_{2}$, and a set $B_{n}\subseteq\mathcal{A}^{n}\times\mathcal{A}^{n}$ such that for every $\mathbf{z}_{n}\in B_{n}$,
\begin{align*}
\mathbb{P}\left(\left.L_{n}\big(\widetilde{\mathbf{Z}}_{n}\big)-L_{n}(\mathbf{Z}_{n})=1\right|\mathbf{Z}_{n}=\mathbf{z}_{n}\right)\geq\varepsilon_{1},\quad\mathbb{P}\left(\left.L_{n}\big(\widetilde{\mathbf{Z}}_{n}\big)-L_{n}(\mathbf{Z}_{n})=-1\right|\mathbf{Z}_{n}=\mathbf{z}_{n}\right)\leq\varepsilon_{2}.
\end{align*}
Moreover, there exists $p_{0}>0$, such that for every $0<p<p_{0}$,
\begin{align*}
\mathbb{P}\left(\mathbf{Z}_{n}\in B_{n}\right)\geq 1-e^{-c_{1}n},
\end{align*}
where $c_{1}>0$ does not depend on $n$, but may depend on $p$.
\end{theorem}
Hence \eqref{eq:ReqRanTrans} holds with $\varepsilon_{0}=\varepsilon_{1}-\varepsilon_{2}$ and $\Delta_{n}(\varepsilon_{0})\leq e^{-c_{1}n}$, for $n$ large enough. Furthermore, as shown above, \eqref{eq:RanTransImpli} also holds for $u_{i}=0,1,\ldots,2n-1$.

In what follows, we will apply the fast convergence result (Theorem \ref{thm:FastConvThm}) to two possible choices of $\mathcal{U}_{n}$, and then compare the corresponding lower bounds on moments and exponential moments of $|L_{n}(\mathbf{Z_{n}})-\mu_{n}|$.

\subsubsection{The Standard $\mathcal{U}_{n}$}\label{subsubsec:StandUnLCS}

Often (cf.~\cite{LemberMatzinger:2009}, \cite{LemberMatzingerTorres:2012(2)} and~\cite{HoudreMa:2016}) the following choice of $\mathcal{U}_{n}$ is used:
\begin{align*}
\mathcal{U}_{n}=\left[2nq-\sqrt{2n},\,2nq+\sqrt{2n}\right]\cap\mathcal{S}_{n}.
\end{align*}
This ``standard $\mathcal{U}_{n}$" is such that $k_{0}=1$ (recall \eqref{eq:DefUnk0}). By the de Moivre-Laplace local limit theorem (cf.~\cite[pp. 56]{Shiryaev:1995}), there exists a universal constant $b=b(p)>0$, which is independent of $n$, but depends on $p$, such that for $n$ large enough,
\begin{align}\label{eq:LowerBoundProbULocalCLT}
\mathbb{P}\left(U_{n}=u\right)\geq\frac{1}{b\sqrt{n}},\quad\text{for any }\,u\in\mathcal{U}_{n}.
\end{align}
The (best) constant $b$ depends on the choice of $\mathcal{U}_{n}$, and it is easy to see (see \eqref{eq:ExtendedUnLowerBound} below) that for the standard $\mathcal{U}_{n}$ as above, $b\geq e^{2}\sqrt{\pi}$. Without loss of generality, assume that $2nq$ is a positive integer (if not, replace $2nq$ in the definition of $\mathcal{U}_{n}$ with $\lceil 2nq\rceil$). Hence, with $\varphi(n)=1/b \sqrt{n}$, the condition \eqref{eq:DeltafastVarphi} is verified and thus, all the assumptions of Theorem \ref{thm:FastConvThm} are satisfied. Therefore, recalling that $\delta=\varepsilon_{0}/2$, $\mathbb{E}(\Phi(|L_{n}(\mathbf{Z}_{n})-\mu_{n}|))$ is lower-bounded by
\begin{align}
\sum_{u_{i}\in\mathcal{U}_{n}}\Phi\left(\delta\left|u_{i}-a_{n}\right|\right)\mathbb{P}\left(U_{n}=u_{i}\right)&\geq\sum_{u_{i}\in\mathcal{U}_{n}}\Phi\left(\delta\left|u_{i}-a_{n}\right|\right)\,\frac{1}{b\sqrt{n}}=\frac{1}{b\sqrt{n}}\sum_{k=2qn-\sqrt{2n}}^{2qn+\sqrt{2n}}\Phi\left(\delta\left|k-a_{n}\right|\right)\nonumber\\
\label{eq:LowerBoundPhi}
&\geq\frac{1}{b\sqrt{n}}\sum_{j=-\sqrt{2n}}^{\sqrt{2n}}\Phi\left(\delta|j|\right)=\frac{1}{b\sqrt{n}}\left(\Phi(0)+2\sum_{j=1}^{\sqrt{2n}}\Phi\left(\frac{\varepsilon_{0}j}{2}\right)\right),
\end{align}
where the last inequality is obtained by minimizing over $a_{n}$ and since $2nq$ is a positive integer.
\begin{example}\label{eg:MomentFastConvStadardU}
For $\Phi(x)=|x|^{r}$, $r\geq 1$, then \eqref{eq:LowerBoundPhi} becomes
\begin{align*}
\frac{\varepsilon_{0}^{r}}{2^{r-1}b\sqrt{n}}\sum_{j=1}^{\sqrt{2n}}j^{r}\geq\frac{\varepsilon_{0}^{r}}{2^{r-1}b\sqrt{n}}\int_{0}^{\sqrt{2n}}x^{r}\,dx=\frac{\varepsilon_{0}^{r}(2n)^{(r+1)/2}}{(r+1)2^{r-1}b\sqrt{n}}=\frac{2^{(3-r)/2}\varepsilon_{0}^{r}}{b\,(r+1)}\,n^{r/2}.
\end{align*}
Hence, for $n$ large enough (independent of $r$, see Remark \ref{rem:UniformLargeN}),
\begin{align}\label{eq:LowerBoundxr}
\mathbb{E}\left(\left|L_{n}(\mathbf{Z}_{n})-\mu_{n}\right|^{r}\right)\geq d_{1}(r)\,n^{r/2},\quad\text{where }\,d_{1}(r)=d_{1}(r,p,\varepsilon_{0}):=\frac{2^{(3-r)/2}\varepsilon_{0}^{r}}{b\,(r+1)}.
\end{align}
\end{example}
\begin{example}\label{eg:ExpMomentFastConvStadardU}
For $\Phi(x)=e^{tx}$, $t>0$, then the second summation in \eqref{eq:LowerBoundPhi} becomes
\begin{align*}
\sum_{j=1}^{\sqrt{2n}}\Phi\left(\frac{\varepsilon_{0}j}{2}\right)=\sum_{j=1}^{\sqrt{2n}}e^{t\varepsilon_{0}j/2}=\sum_{j=1}^{\sqrt{2n}}\rho_{t}^{j}=\frac{\rho_{t}-\rho_{t}^{\sqrt{2n}+1}}{1-\rho_{t}},\quad\text{where }\,\rho_{t}:=e^{\varepsilon_{0}t/2}.
\end{align*}
Hence, for $n$ large enough,
\begin{align}\label{eq:LowerBoundetx}
\mathbb{E}\!\left(e^{\left|L_{n}(\mathbf{Z}_{n})-\mu_{n}\right|t}\right)\!\geq\!\frac{1}{b\sqrt{n}}\!\left(\!1\!+\!\frac{2\rho_{t}\!-\!2\rho_{t}^{\sqrt{2n}+1}}{1-\rho_{t}}\right)\!=\!\frac{\rho_{t}\!\left(2\rho_{t}^{\sqrt{2n}}\!-\!\rho^{-1}_{t}\!-\!1\right)}{b\left(\rho_{t}-1\right)\sqrt{n}}=\frac{2\rho_{t}^{\sqrt{2n}}-1-\rho^{-1}_{t}}{b\left(1-\rho^{-1}_{t}\right)\sqrt{n}}.
\end{align}
Now, for any $t>0$, there exists $n$ large enough (depending on $t$) such that
\begin{align*}
\mathbb{E}\left(e^{\left|L_{n}(\mathbf{Z}_{n})-\mu_{n}\right|t}\right)\geq\frac{2\rho_{t}}{b\left(\rho_{t}-1\right)}\,\rho_{t}^{\sqrt{n}}.
\end{align*}
Taking $t=s/\sqrt{n}$, with $s>0$, gives a lower bound on the moment generating function of $|V_{n}|=|L_{n}(\mathbf{Z}_{n})-\mu_{n}|/\sqrt{n}$. Moreover, by \eqref{eq:LowerBoundetx},
\begin{align}\label{eq:AsyLowerBoundMGFVn}
\liminf_{n\rightarrow\infty}\mathbb{E}\left(e^{s|V_{n}|}\right)\geq\liminf_{n\rightarrow\infty}\frac{2\rho_{t}^{\sqrt{2n}}-1-\rho^{-1}_{t}}{b\left(1-\rho^{-1}_{t}\right)\sqrt{n}}=\frac{4}{b\varepsilon_{0}s}\left(e^{\varepsilon_{0}s/\sqrt{2}}-1\right).
\end{align}
Since
\begin{align*}
\lim_{s\downarrow 0}\frac{1}{\varepsilon_{0}s}\left(e^{\varepsilon_{0}s/\sqrt{2}}-1\right)=\frac{1}{\sqrt{2}},
\end{align*}
the right-hand side of \eqref{eq:AsyLowerBoundMGFVn} approaches $2\sqrt{2}/b$ as $s\downarrow 0$. Since $b\geq\sqrt{\pi}e^{2}$, then $2\sqrt{2}<b$, showing the deficiency of the obtained bound, since we naturally expect the right hand side of \eqref{eq:AsyLowerBoundMGFVn} to approach one as $s\downarrow 0$.
\end{example}
\begin{remark}\label{rem:TaylorApproxMGFVn}
Since $|V_{n}|$ is non-negative, a series expansion, together with lower bounds on moments of $|V_{n}|$, can also lower-bound its moment generating function. Indeed, by \eqref{eq:LowerBoundxr}, for $n$ large enough (independent of $r$),
\begin{align}\label{eq:RefAsyLowerBoundMGFVn}
\mathbb{E}\left(e^{s|V_{n}|}\right)=1+\sum_{r=1}^{\infty}\mathbb{E}\left(|V_{n}|^{r}\right)\frac{s^{r}}{r!}\geq 1+\frac{2\sqrt{2}}{b}\sum_{r=1}^{\infty}\frac{\left(\frac{\varepsilon_{0}s}{\sqrt{2}}\right)^{r}}{(r+1)!}=1-\frac{2\sqrt{2}}{b}+\frac{4}{b\varepsilon_{0}s}\left(e^{\varepsilon_{0}s/\sqrt{2}}-1\right).
\end{align}
This gives a lower bound similar to \eqref{eq:AsyLowerBoundMGFVn}, except for the additive constant $1-2\sqrt{2}/b$, which ensures that the right-hand side of \eqref{eq:RefAsyLowerBoundMGFVn} tends to one as $s\downarrow 0$.

Note finally that \eqref{eq:AsyLowerBoundMGFVn} is computed using the lower bound in \eqref{eq:LowerBoundPhi} with $\Phi(x)=e^{tx}$ and $t=s/\sqrt{n}$:
\begin{align*}
\frac{1}{b\sqrt{n}}\!\sum_{j=-\sqrt{2n}}^{\sqrt{2n}}\!e^{t\varepsilon_{0}|j|/2}=\frac{1}{b\sqrt{n}}\!\sum_{j=-\sqrt{2n}}^{\sqrt{2n}}\sum_{r=0}^{\infty}\frac{1}{r!}\left(\frac{t\varepsilon_{0}|j|}{2}\right)^{r}\!=\frac{2\sqrt{2n}+1}{b\sqrt{n}}+\sum_{r=1}^{\infty}\left[\frac{\varepsilon_{0}^{r}}{2^{r-1}b\sqrt{n}r!}\left(\frac{s}{\sqrt{n}}\right)^{r}\sum_{j=1}^{\sqrt{2n}}j^{r}\right].
\end{align*}
Above, the constant term (i.e., $r=0$) converges to $2\sqrt{2}/b$ as $n\rightarrow\infty$, which is different from the constant term (equal to $1$) in the series expansion in \eqref{eq:RefAsyLowerBoundMGFVn}. The lower bound on the last series above can be computed using \eqref{eq:LowerBoundxr}:
\begin{align*} \sum_{r=1}^{\infty}\left[\frac{\varepsilon_{0}^{r}}{2^{r-1}b\sqrt{n}r!}\left(\frac{s}{\sqrt{n}}\right)^{r}\sum_{j=1}^{\sqrt{2n}}j^{r}\right]\geq\sum_{r=1}^{\infty}\frac{1}{r!}\left(\frac{s}{\sqrt{n}}\right)^{r}\frac{2^{(3-r)/2}\varepsilon_{0}^{r}}{b\,(r+1)}\,n^{\frac{r}{2}}=-\frac{2\sqrt{2}}{b}+\frac{4}{b\varepsilon_{0}s}\left(e^{\varepsilon_{0}s/\sqrt{2}}-1\right).
\end{align*}
Therefore, the difference in the constants in \eqref{eq:AsyLowerBoundMGFVn} and \eqref{eq:RefAsyLowerBoundMGFVn} stems from the different constant terms ($r=0$) in the respective series expansions.
\end{remark}

\subsubsection{The Extended $\mathcal{U}_{n}$}

The bounds presented in the previous subsection give the correct order for the centered absolute moments (given that $p$ is small enough), i.e., $\mathbb{E}(|L_{n}-\mu_{n}|^{r})\asymp O(n^{r/2})$. On the other hand, as argued at the end of Example \ref{eg:ExpMomentFastConvStadardU}, the lower bound on the moment generating function of $|V_{n}|$ can be improved. A way to do so (also valid for the moments of $|V_{n}|$) is to extend the set $\mathcal{U}_{n}$ and to refine the lower bound \eqref{eq:LowerBoundProbULocalCLT}. In the present subsection, we show how the refined lower bound yields better approximations.

In what follows, take $\beta\in(1/2,2/3)$, and redefine $\mathcal{U}_{n}$ as follows
\begin{align}\label{eq:ExtendedUn}
\mathcal{U}_{n}=\left[2nq-(2n)^{\beta},\,2nq+(2n)^{\beta}\right]\cap\mathcal{S}_{n}=\left\{k\in\mathbb{Z}:\,|k-2nq|\leq (2n)^{\beta}\right\}.
\end{align}
In the sequel, this $\mathcal{U}_{n}$ is called the ``extended" $\mathcal{U}_{n}$, and note that again $k_{0}=1$ (recall \eqref{eq:DefUnk0}). Since $(2n)^{\beta}=o(n^{2/3})$, it follows from the local de Moivre-Laplace theorem (cf.~\cite[pp. 56]{Shiryaev:1995}) that
\begin{align*}
\sup_{k:\,|k-2nq|\leq (2n)^{\beta}}\left|\frac{\mathbb{P}\left(U_{n}=k\right)}{\phi_{n}(k)}-1\right|\rightarrow 0,
\end{align*}
where
\begin{align*}
\phi_{n}(k):=\frac{1}{2\sqrt{\pi pqn}}\exp\left(-\frac{\left(k-2nq\right)^{2}}{4npq}\right),\quad\text{and }\,q=1-p.
\end{align*}
Hence, for every $\epsilon>0$, there exists $N\in\mathbb{N}$ such that, for any $n>N$,
\begin{align}\label{eq:LocalDMVIneq}
\left(1+\epsilon\right)\phi_{n}(k)\geq\mathbb{P}\left(U_{n}=k\right)\geq\left(1-\epsilon\right)\phi_{n}(k),\quad\text{for all }\,k\in\mathbb{Z}\,\,\,\,\text{with }\,\left|k-2nq\right|\leq (2n)^{\beta}.
\end{align}
In particular,
\begin{align}\label{eq:ExtendedUnLowerBound}
\mathbb{P}\left(U_{n}=k\right)\geq\frac{\left(1-\epsilon\right)}{2\sqrt{\pi pqn}}\exp\left(-\frac{(2n)^{2\beta-1}}{2pq}\right)=:\varphi(n).
\end{align}
Since $2\beta-1<1$, it follows that
\begin{align*}
\varphi^{-1}(n)\,e^{-c_{1}n}\rightarrow 0,\quad\text{whenever }\,c_{1}>0.
\end{align*}
On the other hand, by Theorem \ref{thm:LCSLargeProb}, there exists $\varepsilon_{0}>0$ and $c_{1}>0$, such that $\Delta_{n}(\varepsilon_{0})\leq e^{-c_{1}n}$, for any $n\in\mathbb{N}$ large enough. Thus, the condition \eqref{eq:DeltafastVarphi} holds and therefore, all the assumptions of Theorem \ref{thm:FastConvThm} are verified.

Thus, recalling that $\delta=\varepsilon_{0}/2$ and $a_{n}\in[2nq-(2n)^{\beta},2nq+(2n)^{\beta}]$,
\begin{align*}
\mathbb{E}\left(\Phi\left(\left|L_{n}(\mathbf{Z}_{n})-\mu_{n}\right|\right)\right)\geq\sum_{k\in\mathcal{U}_{n}}\Phi\left(\delta\left|k-a_{n}\right|\right)\mathbb{P}\left(U_{n}=k\right)\geq(1-\epsilon)\sum_{k:\,|k-2nq|\leq (2n)^{\beta}}\Phi\left(\delta\left|k-a_{n}\right|\right)\phi_{n}(k),
\end{align*}
for $n$ large enough, where the inequality follows from \eqref{eq:LocalDMVIneq}. From the symmetry of $\Phi(\delta|\cdot|)$, the right-hand side above is minimized at $a_{n}=2nq$. Moreover, let
\begin{align*}
j_{k}(n):=\frac{(k-2nq)}{\sqrt{2n}},\quad\text{for }\,k\in\mathbb{Z}\,\text{ with }\,\left|k-2nq\right|\leq (2n)^{\beta}.
\end{align*}
Clearly, $(j_{k}(n))_{k\in\mathcal{U}_{n}}$ forms a $1/\sqrt{2n}$-net over the interval $[-(2n)^{\beta-1/2},(2n)^{\beta-1/2}]$. Hence,
\begin{align}
\mathbb{E}\left(\Phi\left(\left|L_{n}(\mathbf{Z}_{n})-\mu_{n}\right|\right)\right)&\geq (1-\epsilon)\sum_{k:\,|k-2nq|\leq (2n)^{\beta}}\Phi\left(\delta\left|k-2nq\right|\right)\phi_{n}(k)\nonumber\\
\label{eq:LowerBoundPhiExtU} &=(1-\epsilon)\!\sum_{k:\,|k-2nq|\leq (2n)^{\beta}}\frac{1}{2\sqrt{\pi pqn}}\Phi\!\left(\delta\sqrt{2n}\left|j_{k}(n)\right|\right)\exp\left(-\frac{\left(j_{k}(n)\right)^{2}}{2pq}\right).
\end{align}
\begin{example}\label{eg:MomentFastConvExtU}
For $\Phi(x)=|x|^{r}$, $r\geq 1$, since $n^{\beta-1/2}\rightarrow\infty$, as $n\rightarrow\infty$,
\begin{align*}
\lim_{n\rightarrow\infty}\frac{1}{\sqrt{2\pi pq}}\sum_{k:\,|k-2nq|\leq (2n)^{\beta}}\frac{1}{\sqrt{2n}}\left|j_{k}(n)\right|^{r}\exp\left(-\frac{\left(j_{k}(n)\right)^{2}}{2pq}\right)=\frac{1}{\sqrt{2\pi pq}}\int_{-\infty}^{\infty}|x|^{r}e^{-x^{2}/(2pq)}\,dx.
\end{align*}
Since $\epsilon>0$ is arbitrary, \eqref{eq:LowerBoundPhiExtU} leads to
\begin{align*}
\liminf_{n\rightarrow\infty}\mathbb{E}\left(|V_{n}|^{r}\right)=\liminf_{n\rightarrow\infty}\frac{\mathbb{E}\left(\left|L_{n}(\mathbf{Z}_{n})-\mu_{n}\right|^{r}\right)}{n^{r/2}}\geq 2^{r/2}\delta^{r}\mathbb{E}\left(|\xi|^{r}\right)=\frac{\varepsilon_{0}^{r}}{\sqrt{\pi}}(pq)^{r/2}\,\Gamma\left(\frac{r+1}{2}\right)=:d_{2}(r),
\end{align*}
where $\xi\sim\mathcal{N}(0,pq)$. Therefore, for $n\in\mathbb{N}$ large enough,
\begin{align*}
\mathbb{E}\left(\left|L_{n}(\mathbf{Z}_{n})-\mu_{n}\right|^{r}\right)\geq\frac{\varepsilon_{0}^{r}}{\sqrt{\pi}}(pq)^{r/2}\,\Gamma\left(\frac{r+1}{2}\right)n^{r/2}.
\end{align*}

Let us compare this last bound with \eqref{eq:LowerBoundxr} of Example \ref{eg:MomentFastConvStadardU}. From \eqref{eq:LocalDMVIneq},
\begin{align*}
\frac{1}{b}\leq\frac{1}{2\sqrt{\pi pq}}\,e^{-1/(2pq)}\quad\text{implies that}\quad d_{1}(r)\leq\frac{2^{\frac{1-r}{2}}\varepsilon_{0}^{r}}{(r+1)\sqrt{\pi pq}}\,e^{-1/(2pq)}.
\end{align*}
Hence, for each $r\geq 1$, $d_{2}(r)>d_{1}(r)$ will follow from
\begin{align}\label{eq:GammaFunt}
(2pq)^{\frac{r+1}{2}}\Gamma\left(\frac{r+1}{2}\right)>\frac{2}{r+1}\,e^{-1/(2pq)},\quad\text{for }\,p\in(0,1)\,\,\,\text{and}\,\,\,q=1-p.
\end{align}
To show \eqref{eq:GammaFunt}, set $\theta:=(r+1)/2$ and
\begin{align*}
f(x):=x^{-\theta}e^{x}\,\Gamma(\theta+1),\quad x:=\frac{1}{2pq}\geq 2.
\end{align*}
Then \eqref{eq:GammaFunt} is equivalent to $f(x)>1$, for all $x\geq 2$. Since
\begin{align*}
f'(x)=\Gamma(\theta+1)\left(-\theta x^{-\theta-1}e^{x}+x^{-\theta}e^{x}\right)=x^{-\theta}e^{x}\,\Gamma(\theta+1)\left(1-\frac{\theta}{x}\right)
\end{align*}
has a unique zeros at $x=\theta$, it is sufficient to show that
\begin{align*}
f(\theta)=\theta^{-\theta}e^{\theta}\,\Gamma(\theta+1)=\theta^{1-\theta}e^{\theta}\,\Gamma(\theta)>1,\quad\text{for any }\,\theta\geq 1.
\end{align*}
But, this follows from a classical inequality for the gamma function (cf.~\cite[Theorem 1]{KeckicVasic:1971})
\begin{align*}
\frac{\Gamma(x)}{\Gamma(y)}\geq\frac{x^{x-1}e^{y}}{y^{y-1}e^{x}},\quad\text{for any }\,x\geq y\geq 1.
\end{align*}
Therefore, $d_{2}(r)>d_{1}(r)$, for any $r\geq 1$. That is, with the extended $\mathcal{U}_{n}$, Theorem \ref{thm:FastConvThm} provides larger constants than the ones obtained with the standard $\mathcal{U}_{n}$.
\end{example}
\begin{example}\label{eg:ExpMomentFastConvExtU}
For $\Phi_{n}(x)=e^{sx/\sqrt{2n}}$ with $s>0$, then
\begin{align}
\sum_{k:\,|k-2nq|\leq (2n)^{\beta}}\Phi_{n}\left(\delta|k-2nq|\right)\phi_{n}(k)&=\frac{1}{\sqrt{2\pi pq}}\frac{1}{\sqrt{2n}}\sum_{k:\,|k-2nq|\leq (2n)^{\beta}}\exp\left(s\delta\left|j_{k}(n)\right|-\frac{\left(j_{k}(n)\right)^{2}}{2pq}\right)\nonumber\\
\label{eq:LimitExpNormal} &\rightarrow\frac{1}{\sqrt{2\pi pq}}\int_{-\infty}^{\infty}\exp\left(s\delta|x|-\frac{x^{2}}{2pq}\right)dx,\quad n\rightarrow\infty.
\end{align}
Hence, \eqref{eq:LowerBoundPhiExtU} leads to
\begin{align*}
\liminf_{n\rightarrow\infty}\mathbb{E}\left(\exp\left(s\frac{\left|L_{n}(\mathbf{Z}_{n})-\mu_{n}\right|}{\sqrt{2n}}\right)\right)\geq\mathbb{E}\left(e^{s\delta|\xi|}\right),
\end{align*}
where again $\xi\sim\mathcal{N}(0,pq)$. Therefore, for any $s>0$,
\begin{align}\label{eq:LowerBoundVnExtUn}
\liminf_{n\rightarrow\infty}\mathbb{E}\left(e^{s|V_{n}|}\right)\geq\mathbb{E}\left(e^{s\varepsilon_{0}|\xi|/\sqrt{2}}\right).
\end{align}
Unlike the lower bound obtained in \eqref{eq:AsyLowerBoundMGFVn}, the present lower bound approaches one as $s\downarrow 0$. Thus, again, the bound obtained via Theorem \ref{thm:FastConvThm}, with the extended $\mathcal{U}_{n}$, outperforms the one obtained with the standard $\mathcal{U}_{n}$.
\end{example}
\begin{remark}\label{rem:TaylorApproxMGFVnExtUn}
Since $d_{2}$ is given by the $r$-th moment of $|\xi|$, clearly
\begin{align*}
1+\sum_{r=1}^{\infty}d_{2}(r)\frac{s^{r}}{r!}=\mathbb{E}\left(e^{s\varepsilon_{0}|\xi|/\sqrt{2}}\right),
\end{align*}
and so in this case, the lower bound on the moment generating function obtained via the series expansion gives the same result as \eqref{eq:LowerBoundVnExtUn}.
\end{remark}
\begin{remark}
There is another way to get \eqref{eq:LimitExpNormal}. By \eqref{eq:LocalDMVIneq}, it is enough to prove that, as $n\rightarrow\infty$,
\begin{align*}
\sum_{k:\,|k-2nq|\leq
(2n)^{\beta}}\!\!\!\!\Phi_{n}\!\left(\delta\!\left|k-2nq\right|\right)\mathbb{P}\left(U_{n}=k\right)=\!\!\sum_{k:|k-2nq|\leq
(2n)^{\beta}}\!\!\!\!\exp\left\{\frac{s\delta}{\sqrt{2n}}\left|k-2nq\right|\right\}\mathbb{P}\left(U_{n}=k\right)\rightarrow\mathbb{E}\!\left(e^{s\delta|\xi|}\right).
\end{align*}
To see this, note that by the CLT and the continuous mapping theorem, for any $s\in\mathbb{R}$,
\begin{align}\label{eq:WeakConvUn}
\exp\left(\frac{s\delta}{\sqrt{2n}}\left|U_{n}-2nq\right|\right)\convdist e^{s\delta|\xi|},\quad n\rightarrow\infty,
\end{align}
where ``$\convdist$" stands for convergence in law. Moreover, for any $\epsilon>0$, by Hoeffding's exponential inequality,
\begin{align*}
\mathbb{E}\left(\exp\left(\frac{s\delta(1+\epsilon)}{\sqrt{2n}}\left|U_{n}-2nq\right|\right)\right)&=1+\int_{1}^{\infty}\mathbb{P}\left(\exp\left(s\delta(1+\epsilon)\sqrt{2n}\left|\frac{U_{n}-2nq}{2n}\right|\right)\geq x\right)dx\\
&=1+\int_{1}^{\infty}\mathbb{P}\left(\left|\frac{U_{n}}{2n}-q\right|\geq\frac{\ln x}{s\delta(1+\epsilon)\sqrt{2n}}\right)dx\\
&\leq 1+2\int_{1}^{\infty}\exp\left(-4n\frac{(\ln x)^{2}}{2s^{2}\delta^{2}(1+\epsilon)^{2}n}\right)dx\\
&=1+2\int_{0}^{\infty}\exp\left(-\frac{2y^{2}}{s^{2}\delta^{2}(1+\epsilon)^{2}}+y\right)dy<\infty,
\end{align*}
which implies that the family of non-negative random variables
\begin{align*}
\left\{\exp\left(\frac{s\delta}{\sqrt{2n}}\left|U_{n}-2nq\right|\right):\,\,n\in\mathbb{N}\right\}
\end{align*}
is uniformly integrable. Let $A_{n}:=\left\{U_{n}\not\in\mathcal{U}_{n}\right\}$. Again, Hoeffding's exponential inequality leads to
\begin{align}\label{eq:ConvProbAn}
\mathbb{P}\left(A_{n}\right)=\mathbb{P}\left(\left|U_{n}-2nq\right|>(2n)^{\beta}\right)\leq 2\exp\left(-2(2n)^{2\beta-1}\right)\rightarrow 0,\quad
n\rightarrow\infty,
\end{align}
since $2\beta>1$. Therefore, for any $s\in\mathbb{R}$,
\begin{align*}
&\lim_{n\rightarrow\infty}\sum_{k:\,|k-2nq|\leq
(2n)^{\beta}}\exp\left(\frac{s\delta}{\sqrt{2n}}\left|k-2nq\right|\right)\mathbb{P}\left(U_{n}=k\right)\\
&\quad\,=\lim_{n\rightarrow\infty}\mathbb{E}\left({\bf
1}_{A_{n}^{c}}\exp\left(\frac{s\delta}{\sqrt{2n}}\left|U_{n}-2nq\right|\right)\right)=\lim_{n\rightarrow\infty}\mathbb{E}\left(\exp\left(\frac{s\delta}{\sqrt{2n}}\left|U_{n}-2nq\right|\right)\right)=\mathbb{E}\left(e^{s\delta|\xi|}\right),
\end{align*}
where the last equality follows from \eqref{eq:WeakConvUn} and the uniform integrability.
\end{remark}

\subsection{Uniform Approximation}\label{subsec:UnifApproxLCS}

We next show that in the binary LCS case, the uniform approximation approach (Theorem \ref{thm:UnifConv}) also applies provided that $\mathcal{U}_{n}$ is standard. Indeed, with the standard $\mathcal{U}_{n}$, we have $m(n)=2\sqrt{2n}$ (ignoring the rounding), and $\varphi(n)=1/(b\sqrt{n})$, so that the condition $m(n)\geq c\varphi^{-1}(n)=bc\sqrt{n}$ holds with $c=2\sqrt{2}/b$, and thus Theorem \ref{thm:UnifConv} applies. Therefore,
\begin{align}\label{eq:UnifBoundLCS}
\mathbb{E}\left(\Phi\!\left(\left|L_{n}(\mathbf{Z}_{n})\!-\!\mu_{n}\right|\right)\right)\geq\!\!\sum_{j=1}^{\frac{c}{8}\varphi^{-1}(n)}\!\!\Phi\left(\frac{\varepsilon_{0}}{2}\!\left(\frac{c}{8\varphi(n)}\!+\!j\right)\right)\varphi(n)+\frac{c}{8}\Phi\!\left(\frac{\varepsilon_{0}c}{16\varphi(n)}\right)\geq\frac{c}{4}\,\Phi\!\left(\frac{\varepsilon_{0}c}{16\varphi(n)}\right).
\end{align}
Again, let us apply this result to $|L_{n}(\mathbf{Z}_{n})-\mu_{n}|$.
\begin{example}\label{eg:MomentUnifApprox}
For $\Phi(x)=|x|^{r}$, $r\geq 1$, the rightmost bound in \eqref{eq:UnifBoundLCS} gives
\begin{align*}
\mathbb{E}\left(\left|L_{n}(\mathbf{Z}_{n})-\mu_{n}\right|^{r}\right)\geq\frac{c}{4}\left(\frac{\varepsilon_{0}c}{16\varphi(n)}\right)^{r}=\frac{\varepsilon_{0}^{r}}{2^{(5r+1)/2}\,b}\,n^{r/2}.
\end{align*}
Note that
\begin{align*}
d_{3}(r)=d_{3}(r,p,\varepsilon_{0}):=\frac{\varepsilon_{0}^{r}}{2^{(5r+1)/2}\,b}<\frac{2^{\frac{3-r}{2}}\varepsilon_{0}^{r}}{b\,(r+1)}=d_{1}(r),
\end{align*}
where $d_{1}(r)$ is the constant given in \eqref{eq:LowerBoundxr}. This last fact is quite expected since we had a cruder approximation and a weaker assumption $\Delta_{n}(\varepsilon_{0})\rightarrow 0$, as $n\rightarrow\infty$.

Next, applying the finer lower middle bound in \eqref{eq:UnifBoundLCS} gives
\begin{align*}
\sum_{j=1}^{\frac{c}{8}\varphi^{-1}(n)}\!\Phi\!\left(\frac{\varepsilon_{0}}{2}\!\left(\frac{c}{8\varphi(n)}\!+\!j\right)\right)\varphi(n)+\frac{c}{8}\,\Phi\!\left(\frac{\varepsilon_{0}c}{16\varphi(n)}\right)&=\frac{1}{b\sqrt{n}}\left(\frac{\varepsilon_{0}}{2}\right)^{r}\!\sum_{j=1}^{\sqrt{n}/(2\sqrt{2})}\!\left(\frac{\sqrt{n}}{2\sqrt{2}}+j\right)^{r}\!+\frac{\varepsilon_{0}^{r}\,n^{r/2}}{2^{(5r+3)/2}\,b}\\
&\geq\frac{\varepsilon_{0}^{r}}{b\,2^{r}\sqrt{n}}\int_{0}^{\sqrt{n}/(2\sqrt{2})}\!\left(\frac{\sqrt{n}}{2\sqrt{2}}\!+\!x\right)^{r}\!dx+\frac{\varepsilon_{0}^{r}\,n^{r/2}}{2^{(5r+3)/2}\,b}\\
&=\frac{\varepsilon_{0}^{r}}{2^{(5r+3)/2}\,b}\left(\frac{2^{r+1}-1}{r+1}+1\right)n^{r/2}.
\end{align*}
Once more, a better approximation leads to a larger constant
\begin{align*}
d_{4}(r)=d_{4}(r,p,\varepsilon_{0}):=b^{-1}\varepsilon_{0}^{r}\,2^{-\frac{5r+3}{2}}\left(\frac{2^{r+1}-1}{r+1}+1\right)=\frac{2^{r+1}+r}{2(r+1)}\,d_{3}(r)<d_{1}(r).
\end{align*}
\end{example}
\begin{example}\label{eg:ExpMomentUnifApprox}
For $\Phi(x)=e^{tx}$ with $t>0$, the rightmost bound in \eqref{eq:UnifBoundLCS} gives
\begin{align}\label{eq:uniformApproxetx}
\mathbb{E}\left(e^{\left|L_{n}(\mathbf{Z}_{n})-\mu_{n}\right|t}\right)\geq\frac{c}{4}\,\Phi\left(\frac{\varepsilon_{0}c}{16\varphi(n)}\right)=\frac{1}{\sqrt{2}\,b}\,e^{\varepsilon_{0}t\sqrt{n}/(4\sqrt{2})},
\end{align}
while the more refined middle bound, with $\rho_{t}=e^{\varepsilon_{0}t/2}$, gives
\begin{align*}
\sum_{j=1}^{\frac{c}{8}\varphi^{-1}(n)}\!\!\Phi\!\left(\frac{\varepsilon_{0}}{2}\!\left(\frac{c}{8\varphi(n)}\!+\!j\right)\!\right)\!\varphi(n)+\frac{c}{8}\Phi\!\left(\frac{\varepsilon_{0}c}{16\varphi(n)}\right)&=\frac{e^{\varepsilon_{0}t\sqrt{n}/(4\sqrt{2})}}{b\sqrt{n}}\sum_{j=1}^{\sqrt{n}/(2\sqrt{2})}\rho_{t}^{j}+\frac{1}{2\sqrt{2}\,b}\rho_{t}^{\sqrt{n}/(2\sqrt{2})}\\
&=\frac{\rho_{t}}{b\sqrt{n}\left(\rho_{t}\!-\!1\right)}\left(\rho_{t}^{\sqrt{n}/\sqrt{2}}\!-\!\rho_{t}^{\sqrt{n}/(2\sqrt{2})}\right)+\frac{\rho_{t}^{\sqrt{n}/(2\sqrt{2})}}{2\sqrt{2}\,b}.
\end{align*}
For every $t>0$, the convergence to infinity (as $n\rightarrow\infty$) of this last bound is slower than that of the bound \eqref{eq:LowerBoundetx}. Hence, again, the uniform bound is smaller (for large $n$) than the ones obtained in the previous subsection using Theorem \ref{thm:FastConvThm}. Since $t>0$, for $n$ large enough (depending on $t$),
\begin{align*}
\mathbb{E}\left(e^{\left|L_{n}(\mathbf{Z}_{n})-\mu_{n}\right|t}\right)\geq\frac{1}{b}\,e^{\varepsilon_{0}t\sqrt{n}/(4\sqrt{2})}\,\frac{\rho_{t}^{\frac{1}{4}\sqrt{n}+1}}{\rho_{t}-1}+\frac{e^{\varepsilon_{0}t\sqrt{n}/(4\sqrt{2})}}{2\sqrt{2}\,b}\geq\frac{1}{b}\,e^{\varepsilon_{0}t\sqrt{n}/4}+\frac{e^{\varepsilon_{0}t\sqrt{n}/(4\sqrt{2})}}{2\sqrt{2}\,b}.
\end{align*}
\end{example}
\begin{remark}\label{remark:TaylorUniformApprox}
Again, when estimating the moment generating function of $|V_{n}|=|L_{n}(\mathbf{Z})-\mu_{n}|/\sqrt{n}$ via its series expansion, with the lower bound $\mathbb{E}(|V_{n}|^{r})\geq d_{3}(r)$, for large $n$ (uniformly in $r$),
\begin{align*}
\mathbb{E}\!\left(e^{s|V_{n}|}\right)\geq 1+\frac{1}{\sqrt{2}b}\sum_{r=1}^{\infty}\frac{1}{r!}\left(\frac{\varepsilon_{0}}{4\sqrt{2}}s\right)^{r}=\frac{1}{\sqrt{2}\,b}\,e^{\varepsilon_{0}s/(4\sqrt{2})}+1-\frac{1}{\sqrt{2}b},
\end{align*}
which is the same bound as in \eqref{eq:uniformApproxetx} with $t=s/\sqrt{n}$, except for the constant term.
\end{remark}
\begin{remark}\label{rem:NoUnifApproxExtUn}
With the extended $\mathcal{U}_{n}$, the uniform approximation (Theorem \ref{thm:UnifConv}) cannot be applied, since there does not exist a constant $c>0$ such that $m(n)=2(2n)^{\beta}\geq c\,\varphi^{-1}(n)$, for any $n\in\mathbb{N}$, with $\varphi(n)$ given in \eqref{eq:ExtendedUnLowerBound}.
\end{remark}
\begin{remark}\label{rem:EstMoment}
We have obtained several constants $d_{3}(r)<d_{4}(r)<d_{1}(r)<d_{2}(r)$ involved in lower bounding the moment $\mathbb{E}(|V_{n}|^{r})$. The best constant, $d_{2}(r)$, is obtained under fast convergence with extended $\mathcal{U}_{n}$. Clearly, in order to deduce the existence of $c(r)$ such that $\mathbb{E}(|V_{n}|^{r})\geq c(r)$, it suffices to tackle the case $r=1$. However, it is easy to verify that in all the cases considered above (i.e., $c(r)$ being either $d_{1}(r)$, $d_{2}(r)$, $d_{3}(r)$ or $d_{4}(r)$), the constant $(c(1))^{r}$ is smaller than $c(r)$. Moreover, in all these cases, it is also true that $c(r+1)\geq (c(r))^{(r+1)/r}$, so that a better constant is obtained when estimating the $r$-th moment directly rather than estimating a lower-order moment.
\end{remark}

\section{An Upper Bound on the Rate Function}\label{sec:Rate}

\subsection{Background and Preliminary Results}

The analysis right below (before Proposition \ref{prop:LowerBoundMGFLnl2nq}) is similar to~\cite[Theorem 2]{ArratiaWaterman:1994}, but with a more general notion of score function and a slightly different definition of gap\footnote{In~\cite{ArratiaWaterman:1994}, the scoring function takes the value $1$ for matches and the penalty $-\mu$ for mismatches. Moreover, the gap in~\cite{ArratiaWaterman:1994} is an indel in one sequence, and the gap price $-\delta$ is assumed to be negative.}. Again, let $(X_{n})_{n\in\mathbb{N}}$ and $(Y_{n})_{n\in\mathbb{N}}$ be two independent sequences of i.i.d. random variables, and let us consider a general scoring function $S:\mathcal{A}\times\mathcal{A}\rightarrow\mathbb{R}^{+}$. By subadditivity,
\begin{align*}
L_{n+m}\geq L_{n}+L\left(\mathbf{X}_{n+1}^{n+m};\mathbf{Y}_{n+1}^{n+m}\right),
\end{align*}
where $L(\mathbf{X}_{n+1}^{n+m};\mathbf{Y}_{n+1}^{n+m})$ denotes the optimal alignments score of $\mathbf{X}_{n+1}^{n+m}:=(X_{n+1},\ldots,X_{n+m})$ and $\mathbf{Y}_{n+1}^{n+m}:=(Y_{n+1},\ldots,Y_{n+m})$. Hence, for any $s\geq 0$,
\begin{align*}
\mathbb{P}\left(L_{n+m}\geq s(n+m)\right)\geq\mathbb{P}\left(L_{n}\geq sn,\,L\left(\mathbf{X}_{n+1}^{n+m};\mathbf{Y}_{n+1}^{n+m}\right)\geq sm\right)=\mathbb{P}\left(L_{n}\geq sn\right)\mathbb{P}\left(L_{m}\geq sm\right).
\end{align*}
Thus, by Fekete's lemma, for any $s\geq 0$, the following limit -- the {\it rate function} -- exists:
\begin{align*}
r(s):=\lim_{n\rightarrow\infty}-\frac{1}{n}\ln\mathbb{P}\left(L_{n}\geq sn\right)=\inf_{n\in\mathbb{N}}-\frac{1}{n}\ln\mathbb{P}\left(L_{n}\geq sn\right)<\infty.
\end{align*}
and so,
\begin{align*}
\mathbb{P}\left(L_{n}\geq sn\right)\leq e^{-r(s)\,n},\quad\text{for any }\,s\geq 0,\,\,n\in\mathbb{N}.
\end{align*}
Since for any $n\in\mathbb{N}$,
\begin{align*}
\mathbb{E}\left(L_{n}\right)\leq\gamma^{*}n,\quad\text{where }\,\gamma^{*}:=\lim_{n\rightarrow\infty}\frac{\mathbb{E}\left(L_{n}\right)}{n},
\end{align*}
it follows from \eqref{eq:HoeffdingIneqOptScore} that for any $s>\gamma^{*}$,
\begin{align*}
\mathbb{P}\left(L_{n}\geq
sn\right)=\mathbb{P}\!\left(\!L_{n}\!-\!\mathbb{E}(L_{n})\!\geq\!\left(s\!-\!\frac{\mathbb{E}(L_{n})}{n}\right)\!n\!\right)\leq\mathbb{P}\left(L_{n}\!-\!\mathbb{E}(L_{n})\!\geq\!\left(s\!-\!\gamma^{*}\right)n\right)\leq\exp\!\left(\!-\frac{\left(s-\gamma^{*}\right)^{2}n}{K^{2}}\!\right).
\end{align*}
Therefore, for $s>\gamma^{*}$,
\begin{align}\label{eq:LowerBoundRateFunt}
r(s)\geq\frac{\left(s-\gamma^{*}\right)^{2}}{K^{2}}.
\end{align}
Moreover, for $s=\gamma^{*}$,
\begin{align*}
\mathbb{P}\left(L_{n}\geq\gamma^{*}n\right)=\mathbb{P}\left(L_{n}-\mathbb{E}(L_{n})\geq\left(\gamma^{*}-\frac{\mathbb{E}(L_{n})}{n}\right)n\right)\leq\exp\left(-\frac{n}{K^{2}}\left(\gamma^{*}-\frac{\mathbb{E}(L_{n})}{n}\right)^{2}\right).
\end{align*}
The aim of the present section is to show that the methodology developed to date allows us to partially reverse \eqref{eq:LowerBoundRateFunt} in the LCS case. That is, we will show that there exists a universal constant $B>0$, such that
\begin{align}\label{eq:ConjCond}
r(s)\leq B(s-\gamma^{*})^{2},
\end{align}
for any $s$ belonging to an upper neighborhood of $\gamma^{*}$. Moreover, the full claim \eqref{eq:ConjCond}, i.e., for all $s>\gamma^{*}$, is shown to hold under a uniform convergence assumption, which is not proved in general.

Hereafter and throughout this section, we will only consider the binary LCS case. For the rest of this subsection, we obtain a lower bound on $\mathbb{E}(e^{t(L_{n}(\mathbf{Z}_{n})-\ell_{n}(2nq))})$ (note that there is no absolute value in the exponent), which will be used to prove the claim \eqref{eq:ConjCond} in the neighborhood of $\gamma^{*}$. Since we are considering the special case of LCS and since we have already seen the advantage of the extended $\mathcal{U}_{n}$ over the standard one, we shall use the extended $\mathcal{U}_{n}$. Throughout this section, let $p_{0}>0$ be given as in Theorem \ref{thm:LCSLargeProb}.
\begin{proposition}\label{prop:LowerBoundMGFLnl2nq}
Let $t>0$ and let $p\in(0,p_{0})$. Then, there exists $\varepsilon_{0}>0$, such that
\begin{align*}
\liminf_{n\rightarrow\infty}\mathbb{E}\left(\exp\left(\frac{t}{\sqrt{2n}}\left(L_{n}(\mathbf{Z}_{n})-\ell_{n}(2nq)\right)\right)\right)\geq\lambda\,e^{\,pq\,\varepsilon_{0}^{2}t^{2}/8},
\end{align*}
where
\begin{align}\label{eq:ConstA}
\lambda=\lambda(\varepsilon_{0}):=\min_{t>0}\left(1-\,\mathbb{P}\left(\frac{pq\,\varepsilon_{0}t}{2}<\xi<pqt\right)\right)\in(0,1),
\end{align}
and where $\xi\sim\mathcal{N}(0,pq)$.
\end{proposition}

\noindent
\textbf{Proof.} For simplicity, assume $2nq$ to be an integer. Let $\mathcal{U}_{n}$ be given by \eqref{eq:ExtendedUn}. By the conditional Jensen's inequality,
\begin{align}
\mathbb{E}\!\left(e^{t\left(L_{n}(\mathbf{Z}_{n})-\ell_{n}(2nq)\right)}\right)&\geq\mathbb{E}\left(e^{t\left(\ell_{n}(U_{n})-\ell_{n}(2nq)\right)}\right)\geq\sum_{k\in\mathcal{U}_{n}}e^{t\left(\ell_{n}(k)-\ell_{n}(2nq)\right)}\mathbb{P}\left(U_{n}=k\right)\nonumber\\
\label{eq:DecomMGFLZl2pn} &=\!\sum_{k=2nq}^{2nq+(2n)^{\beta}}\!\!\!e^{t\left(\ell_{n}(k)-\ell_{n}(2nq)\right)}\mathbb{P}\!\left(U_{n}\!=\!k\right)\!+\!\!\!\sum_{k=2nq-(2n)^{\beta}}^{2nq-1}\!\!\!\!\!e^{t\left(\ell_{n}(k)-\ell_{n}(2nq)\right)}\mathbb{P}\!\left(U_{n}\!=\!k\right).
\end{align}
Theorem \ref{thm:LCSLargeProb} ensures that, with $\varepsilon_{0}:=\varepsilon_{1}-\varepsilon_{2}$, the condition \eqref{eq:DeltafastVarphi} holds and therefore, all the assumptions of Theorem \ref{thm:FastConvThm} are satisfied. If $k\in\mathcal{U}_{n}$ and $k\geq 2nq$, then by Theorem \ref{thm:FastConvThm}, $\ell_{n}(k)-\ell_{n}(2nq)\geq\delta(k-2nq)$, where $\delta=\varepsilon_{0}/2<1$; while if $k<2nq$, then by \eqref{eq:star}, $\ell_{n}(k)-\ell_{n}(2nq)\geq k-2nq$. Hence, by \eqref{eq:LocalDMVIneq}, for any $\epsilon>0$ and $n$ (depending on $\epsilon$, but not on $t$) large enough,
\begin{align*}
\sum_{k=2nq}^{2nq+(2n)^{\beta}}\exp\left(\frac{t}{\sqrt{2n}}\left(\ell_{n}(k)-\ell_{n}(2nq)\right)\right)\mathbb{P}\left(U_{n}=k\right)\geq (1-\epsilon)\sum_{k=2nq}^{2nq+(2n)^{\beta}}\exp\left(\frac{t\delta}{\sqrt{2n}}\left(k-2nq\right)\right)\phi_{n}(k),
\end{align*}
and
\begin{align*}
\sum_{k=2nq}^{2nq+(2n)^{\beta}}\exp\left(\frac{t\delta}{\sqrt{2n}}\left(k-2nq\right)\right)\phi_{n}(k)\rightarrow\frac{1}{2}\left(1+\text{erf}\left(\frac{\sqrt{pq}\,\delta t}{\sqrt{2}}\right)\right)e^{\,pq\,\varepsilon_{0}^{2}t^{2}/8}.
\end{align*}
Similarly, for $n$ (depending on $\epsilon$, but not on $t$) large enough,
\begin{align*}
\sum_{k=2nq-(2n)^{\beta}}^{2nq-1}\!\!\!\exp\left(\frac{t}{\sqrt{2n}}\left(\ell_{n}(k)-\ell_{n}(2nq)\right)\right)\mathbb{P}\left(U_{n}=k\right)\geq (1-\epsilon)\!\!\!\sum_{k=2nq-(2n)^{\beta}}^{2nq-1}\!\!\!\exp\left(\frac{t}{\sqrt{2n}}\left(k-2nq\right)\right)\phi_{n}(k),
\end{align*}
and
\begin{align*}
\sum_{k=2nq-(2n)^{\beta}}^{2nq-1}\exp\left(\frac{t}{\sqrt{2n}}\left(k-2nq\right)\right)\phi_{n}(k)\rightarrow\frac{1}{2}\left(1-\text{erf}\left(\frac{\sqrt{pq}\,t}{\sqrt{2}}\right)\right)e^{pq\,t^{2}/2}.
\end{align*}
Since $\epsilon>0$ is arbitrarily chosen,
\begin{align*}
&\liminf_{n\rightarrow\infty}\,\mathbb{E}\left(\exp\left(\frac{t}{\sqrt{2n}}\left(L_{n}(\mathbf{Z}_{n})-\ell_{n}(2nq)\right)\right)\right)\\
&\quad\,\geq\frac{1}{2}\left(1+\text{erf}\left(\frac{\sqrt{pq}\,\delta t}{\sqrt{2}}\right)\right)e^{pq\,\varepsilon_{0}^{2}t^{2}/8}+\frac{1}{2}\left(1-\text{erf}\left(\frac{\sqrt{pq}\,t}{\sqrt{2}}\right)\right)e^{pq\,t^{2}/2}\\
&\quad\,>\left[1+\frac{1}{2}\left(\text{erf}\left(\frac{\sqrt{pq}\,\delta t}{\sqrt{2}}\right)-\text{erf}\left(\frac{\sqrt{pq}\,t}{\sqrt{2}}\right)\right)\right]e^{pq\,\varepsilon_{0}^{2}t^{2}/8}\\
&\quad\,=\left(1-\mathbb{P}\left(pq\,\delta t<\xi<pq\,t\right)\right)e^{pq\,\varepsilon_{0}^{2}t^{2}/8},
\end{align*}
where $\xi\sim N(0,pq)$.\hfill $\Box$

\newpage
To conclude this subsection, we provide an estimate on $|\mathbb{E}(L_{n}(\mathbf{Z}_{n}))-\ell_{n}(2nq)|/\sqrt{2n}$, which will be useful in the sequel.
\begin{corollary}\label{cor:LowerBoundLnZnln2nq}
In the binary LCS setting, with $p\in(0,p_{0})$,
\begin{align*}
\limsup_{n\rightarrow\infty}\frac{\left|\mathbb{E}\left(L_{n}(\mathbf{Z}_{n})\right)-\ell_{n}(2nq)\right|}{\sqrt{2n}}\leq\sqrt{\frac{2pq}{\pi}}.
\end{align*}
In particular,
\begin{align}\label{eq:LimitLn}
\lim_{n\rightarrow\infty}\frac{\ell_{n}(2nq)}{n}=\lim_{n\rightarrow\infty}\frac{\mathbb{E}\left(L_{n}\right)}{n}=\gamma^{*}.
\end{align}
\end{corollary}

\noindent
\textbf{Proof.} By \eqref{eq:star},
\begin{align*}
\left|\ell_{n}(k)-\ell_{n}(2nq)\right|\leq\left|k-2nq\right|,\quad k\in\mathcal{U}_{n}.
\end{align*}
Hence, for any $\epsilon>0$ and $n$ large enough,
\begin{align*}
\sum_{k\in\mathcal{U}_{n}}\frac{\left|\ell_{n}(k)-\ell_{n}(2nq)\right|}{\sqrt{2n}}\mathbb{P}\left(U_{n}=k\right)\leq
(1+\epsilon)\sum_{k\in\mathcal{U}_{n}}\frac{\left|k-2nq\right|}{\sqrt{2n}}\phi_{n}(k).
\end{align*}
Again,
\begin{align*}
\sum_{k\in\mathcal{U}_{n}}\frac{\left|k-2nq\right|}{\sqrt{2n}}\phi_{n}(k)\rightarrow\int_{-\infty}^{\infty}\frac{|x|}{\sqrt{2\pi pq}}\,e^{-x^{2}/(2pq)}\,dx=\sqrt{\frac{2pq}{\pi}},
\end{align*}
and by \eqref{eq:ConvProbAn},
\begin{align*}
\sum_{k\not\in\mathcal{U}_{n}}\frac{\left|\ell_{n}(k)-\ell_{n}(2nq)\right|}{\sqrt{2n}}\mathbb{P}\left(U_{n}=k\right)\leq\frac{n}{\sqrt{2n}}\mathbb{P}\left(U_{n}\not\in\mathcal{U}_{n}\right)\leq\sqrt{2n}\exp\left[-2(2n)^{2\beta-1}\right]\rightarrow 0,\quad n\rightarrow\infty.
\end{align*}
Thus,
\begin{align*}
\limsup_{n\rightarrow\infty}\frac{\left|\mathbb{E}\left(L_{n}(\mathbf{Z}_{n})\right)-\ell_{n}(2nq)\right|}{\sqrt{2n}}\leq\limsup_{n\rightarrow\infty}\frac{\mathbb{E}\left(\left|\ell_{n}(U_{n})-\ell_{n}(2nq)\right|\right)}{\sqrt{2n}}\leq\sqrt{\frac{2pq}{\pi}},
\end{align*}
which clearly implies \eqref{eq:LimitLn}.\hfill $\Box$

\subsection{An Upper Bound on the Rate Function in the Neighborhood of $\gamma^{*}$}\label{subsec:UpperBoundRateFunt}

The goal of this subsection is to prove the claim \eqref{eq:ConjCond} in the neighborhood of $\gamma^{*}$, i.e., for some constant $B>0$ and $\widetilde{\gamma}>0$,
\begin{align}\label{eq:ConjCondLocal}
r(s)\leq B\left(s-\gamma^{*}\right)^{2},\quad\text{for any }\,s\in(\gamma^{*},\gamma^{*}+\widetilde{\gamma}].
\end{align}
Let
\begin{align*}
\Lambda_{n}(t):=\log\mathbb{E}\left(e^{tL_{n}(\mathbf{Z}_{n})}\right),\quad t\in\mathbb{R}.
\end{align*}
The following result of~\cite{Hammersley:1974} (see also~\cite[Theorem 1]{GrossmannYakir:2004}) will be useful in the sequel of the proof: there exists $\Lambda(t)$, $t\in\mathbb{R}$, such that
\begin{align*}
\lim_{n\rightarrow\infty}\frac{\Lambda_{n}(t)}{n}=\Lambda(t),\quad\text{for any }\,t\in\mathbb{R}.
\end{align*}
Moreover, $r(s)$ and $\Lambda(t)$ are convex functions and are
related via
\begin{align}\label{eq:RelatRateLambda}
r(s)=\Lambda^{*}(s):=\sup_{t\geq
0}\left(ts-\Lambda(t)\right),\quad\Lambda(t)=r^{*}(t)=\sup_{s\in\mathbb{R}}\left(ts-r(s)\right),\quad\text{for all }s\in\mathbb{R}\,\,\,\text{and}\,\,\,t\geq 0.
\end{align}

\subsubsection{A Global Upper Bound Under a Uniform Lower Bound}

Recall from Proposition \ref{prop:LowerBoundMGFLnl2nq} that
\begin{align*}
\liminf_{n\rightarrow\infty}\,\mathbb{E}\left[\exp\left(\frac{t}{\sqrt{2n}}\left(L_{n}(\mathbf{Z})-\ell_{n}(2nq)\right)\right)\right]\geq\lambda\,e^{pq\,\varepsilon_{0}^{2}t^{2}/8},
\end{align*}
where $\lambda$ is given by \eqref{eq:ConstA}. Hence, for every $0<\lambda_{1}<\lambda$ and $t>0$,
\begin{align}\label{eq:LowerBounda1}
\mathbb{E}\left[\exp\left(\frac{t}{\sqrt{2n}}\left(L_{n}(\mathbf{Z}_{n})-\ell_{n}(2nq)\right)\right)\right]\geq\lambda_{1}\,e^{pq\,\varepsilon_{0}^{2}t^{2}/8},
\end{align}
provided that $n$ (which depends on $t$) is large enough .
\begin{proposition}\label{lem:GlobUpBoundRateFunt}
Let \eqref{eq:LowerBounda1} hold uniformly over $[0,\infty)$, i.e., let there exists $n_{0}\in\mathbb{N}$ (independent of $t$) such
that, whenever $n\geq n_{0}$, \eqref{eq:LowerBounda1} holds true for every $t\geq 0$, then \eqref{eq:ConjCond} holds true for all $s>\gamma^{*}$.
\end{proposition}
\noindent
\textbf{Proof.} Let $\lambda_{2}:=pq\,\varepsilon_{0}^{2}/8$. Since \eqref{eq:LowerBounda1} holds uniformly in $t$, for $n\geq n_{0}$,
\begin{align*}
\Lambda_{n}\left(\frac{t}{\sqrt{2n}}\right)\geq\frac{\ell_{n}(2nq)t}{\sqrt{2n}}+\ln\lambda_{1}+\lambda_{2}\,t^{2},\quad\text{for any }\,t\geq 0.
\end{align*}
With $u=t/\sqrt{2n}$, for $n\geq n_{0}$,
\begin{align*}
\frac{\Lambda_{n}(u)}{n}\geq\frac{\ell_{n}(2nq)u}{n}+\frac{\ln\lambda_{1}}{n}+2\lambda_{2}u^{2},\quad\text{for any }\,u\geq 0,
\end{align*}
and thus (recalling \eqref{eq:LimitLn}),
\begin{align*}
\Lambda(u)=\lim_{n\rightarrow\infty}\frac{\Lambda_{n}(u)}{n}\geq\gamma^{*}u+2\lambda_{2}u^{2},\quad\text{for any }\,u\geq 0.
\end{align*}
Now, for every $s>\gamma^{*}$, $u\geq 0$,
\begin{align*}
su-\Lambda(u)\leq\left(s-\gamma^{*}\right)u-2\lambda_{2}u^{2},
\end{align*}
which implies that,
\begin{align*}
r(s)=\sup_{u\geq 0}\left(su-\Lambda(u)\right)\leq\sup_{u\geq 0}\left[\left(s-\gamma^{*}\right)u-2\lambda_{2}u^{2}\right]=\frac{\left(s-\gamma^{*}\right)^{2}}{8\lambda_{2}}.
\end{align*}
Therefore, \eqref{eq:ConjCond} holds true, for all $s>\gamma^{*}$, with $B=1/(8\lambda_{2})=1/(4pq\delta^{2})$.\hfill $\Box$

\subsubsection{Binomial Approximation With the Further Extended $\mathcal{U}_{n}$}

In this part, we derive \eqref{eq:ConjCondLocal} without assuming that \eqref{eq:LowerBounda1} holds uniformly in $t$. Let $M_{q}(t)$ and $K_{q}(t)$ be respectively the moment generating function and the cumulant generating function of the Rademacher law with parameter $q$, i.e.,
\begin{align*}
M_{q}(t):=e^{t(1-q)}q+e^{-tq}(1-q),\quad K_{q}(t):=\ln M_{q}(t)=\ln\left(e^{t(1-q)}q+e^{-tq}(1-q)\right).
\end{align*}
We start with the following general theorem.
\begin{theorem}\label{thm:LocUpBoundRateFunt}
Let $p\in(0,p_{0})$. Let there exist $t_{0}>0$ and $\delta>0$ such that, for any $t\in[0,t_{0}]$, whenever $n$ is large enough (possibly depending on $t$),
\begin{align}\label{eq:lambda}
\mathbb{E}\left(e^{t\left(L_{n}(\mathbf{Z}_{n})-\ell_{n}(2nq)\right)}\right)\geq\nu_{n}M^{2n}_{q}(\delta t),
\end{align}
where $(\nu_{n})_{n\in\mathbb{N}}\subseteq [0,\infty)$ is such that $\ln\nu_{n}=o(n)$, as $n\rightarrow\infty$. Then, there exist constants $B>0$ and $\widetilde{\gamma}>0$ such that \eqref{eq:ConjCondLocal} holds true.
\end{theorem}
\noindent
\textbf{Proof.} Let $t\in [0,t_{0}]$. Taking logarithms in \eqref{eq:lambda} leads to
\begin{align*}
\Lambda_{n}(t)=\ln\mathbb{E}\left(e^{tL_{n}(\mathbf{Z}_{n})}\right)\geq t\,\ell_{n}(2nq)+2n\,K_{q}(\delta t).
\end{align*}
Dividing both sides of the above equality by $n$, and using \eqref{eq:LimitLn}, lead to
\begin{align}\label{lambdabound}
\Lambda(t)\geq t\gamma^{*}+2K_{q}(\delta t).
\end{align}
Now for any $s\in\mathbb{R}$,
\begin{align*}
ts-\Lambda(t)\leq ts-t\gamma^{*}-2K_{q}(\delta t)=t(s-\gamma^{*})-2K_{q}(\delta t),
\end{align*}
which implies that
\begin{align}\label{eq:SupLambda0t0}
\sup_{t\in [0,t_{0})}\left(ts-\Lambda(t)\right)\leq\sup_{t\in [0,t_{0})}\left[t\left(s-\gamma^{*}\right)-2K_{q}(\delta t)\right].
\end{align}

We next show that $\Lambda'(0+)=\gamma^{*}$ (recall that $\Lambda(t)$ and $r(s)$ are related to each other only for $t\in[0,\infty)$, see \eqref{eq:RelatRateLambda}). First, by \eqref{lambdabound}, since $K_{q}(\delta t)\geq 0$ for any $t\geq 0$, we have
\begin{align}\label{eq:LowerBoundLambdat}
\Lambda(t)\geq\gamma^{*}t,\quad\text{for any }\,t\geq 0.
\end{align}
Next, by definition, $r(s)=0$ for any $s<\gamma^{*}$. As stated in~\cite[Theorem 1]{GrossmannYakir:2004}, $r$ is a convex and, therefore, continuous function, and so $r(\gamma^{*})=0$. Thus, the function $r$ is equal to zero up to $\gamma^{*}$ and is strictly increasing afterwards, and so,
\begin{align*}
\Lambda(t)=\sup_{s\in\mathbb{R}}\left(ts-r(s)\right)=\sup_{s>\gamma^{*}}\left(ts-r(s)\right).
\end{align*}
Using \eqref{eq:LowerBoundRateFunt} with $K=1$ in the LCS case, for any $t>0$,
\begin{align}\label{eq:UpperBoundLambdat}
\Lambda(t)=\sup_{s>\gamma^{*}}\left(ts-r(s)\right)\leq\sup_{s>\gamma^{*}}\left[ts-\left(s-\gamma^{*}\right)^{2}\right]=ts_{t}-\left(s_{t}-\gamma^{*}\right)^{2}=t\gamma^{*}+\frac{1}{4}t^{2},
\end{align}
where $s_{t}:=t/2+\gamma^{*}>\gamma^{*}$, $t>0$. Combining \eqref{eq:LowerBoundLambdat} and \eqref{eq:UpperBoundLambdat}, it follows that $\Lambda'(0+)=\gamma^{*}$, which, together with \eqref{eq:SupLambda0t0}, implies that there exists $\widetilde{\gamma}_{1}>0$, such that for any $s\in(\gamma^{*},\gamma^{*}+\widetilde{\gamma}_{1}]$,
\begin{align*}
r(s)=\sup_{t\geq 0}\left(ts-\Lambda(t)\right)=\sup_{t\in[0,t_{0}]}\left(ts-\Lambda(t)\right)\leq\sup_{t\in[0,t_{0}]}\left[t\left(s-\gamma^{*}\right)-2K_{q}(\delta t)\right].
\end{align*}
Now,
\begin{align*}
\sup_{t\in [0,t_{0}]}\left[t\left(s-\gamma^{*}\right)-2K_{q}(\delta t)\right]&=2\sup_{t\in [0,t_{0}]}\left[\delta t\left(q+\frac{s-\gamma^{*}}{2\delta}\right)-\ln\left(qe^{t\delta}+(1-q)\right)\right]\\
&=2\sup_{u\in[0,\delta t_{0}]}\left[u\left(q+\frac{s-\gamma^{*}}{2\delta}\right)-\ln\left(qe^{u}+(1-q)\right)\right].
\end{align*}
With
\begin{align*}
x:=q+\frac{s-\gamma^{*}}{2\delta},
\end{align*}
we obtain that the solution of
\begin{align*}
x=\frac{qe^{u}}{qe^{u}+(1-q)}
\end{align*}
is
\begin{align*}
u(x)=\ln\frac{x(1-q)}{q(1-x)}.
\end{align*}
Clearly, there exists $\widetilde{\gamma}_{2}>0$, such that when $s-\gamma^{*}\leq\widetilde{\gamma}_{2}$, $u(x)\in [0,\delta t_{0}]$, and in this case
\begin{align*}
\sup_{u\in[0,\delta t_{0})}\left[u\left(q+\frac{s-\gamma^{*}}{2\delta}\right)-\ln\left(qe^{u}+(1-q)\right)\right]&=\sup_{u\geq 0}\left[u\left(q+\frac{s-\gamma^{*}}{2\delta}\right)-\ln\left(qe^{u}+(1-q)\right)\right]\\
&=u(x)x-\ln\left(qe^{u(x)}+(1-q)\right)\\&=x\ln{x\over q}+(1-x)\ln {1-x\over 1-q}\\
&=:D(x||q).
\end{align*}
Since for $w>0$, $\ln(1+w)<w$, and since for $w\in(0,1)$, $\ln(1-w)<-w$, it follows that
\begin{align*}
D(x||q)\leq x\,\frac{x-q}{q}-(1-x)\frac{x-q}{1-q}=\frac{(x-q)^{2}}{q(1-q)}=\frac{1}{4\delta^{2}q(1-q)}\left(s-\gamma^{*}\right)^{2}.
\end{align*}
Therefore, \eqref{eq:ConjCondLocal} holds with $B=1/(4\delta^{2}q(1-q))$ and $\widetilde{\gamma}=\min(\widetilde{\gamma}_{1},\widetilde{\gamma}_{2})$.\hfill $\Box$
\begin{remark}\label{rem:GlobLowerBoundRateFunt}
By examining the above proof, it is easy to see that if \eqref{eq:lambda} holds for any $t>0$, then \eqref{eq:ConjCond} holds for any $s>\gamma^{*}$.
\end{remark}
Let us return to \eqref{eq:lambda} and examine under which conditions the hypothesis of Theorem \ref{thm:LocUpBoundRateFunt} are satisfied. Since $U_{n}\sim B(2n,q)$,
\begin{align*}
\sum_{k=0}^{2n}e^{t\delta\left(k-2nq\right)}\mathbb{P}\left(U_{n}=k\right)=M_{q}^{2n}(\delta t),
\end{align*}
one might try to use the first term in \eqref{eq:DecomMGFLZl2pn} to get \eqref{eq:lambda}:
\begin{align*}
\mathbb{E}\left(e^{t\left(L_{n}(\mathbf{Z}_{n})-\ell_{n}(2nq)\right)}\right)\geq\sum_{k\in\mathcal{U}_{n},\,k\geq 2nq}e^{t\left(\ell_{n}(k)-\ell_{n}(2nq)\right)}\mathbb{P}\left(U_{n}=k\right).
\end{align*}
Thus \eqref{eq:lambda} holds if,
\begin{align*}
\sum_{k=2nq}^{2nq+(2n)^{\beta}}e^{t\delta\left(k-2nq\right)}\mathbb{P}\left(U_{n}=k\right)\geq\nu_{n}\,M_{q}^{2n}(\delta t),
\end{align*}
for a sequence $(\nu_{n})_{n\in\mathbb{N}}$ of positive reals such that $\ln\nu_{n}=o(n)$ (as $n\rightarrow\infty$). Unfortunately, there does not exist such $(\nu_{n})_{n\in\mathbb{N}}$, since
\begin{align*}
\frac{1}{n}\ln\left(\sum_{k=2nq}^{2nq+(2n)^{\beta}}e^{t\delta(k-2nq)}\mathbb{P}\left(U_{n}=k\right)\right)<\frac{1}{n}\ln\left(\exp\left(t\delta (2n)^{\beta}\right)\right)\rightarrow 0,\quad n\rightarrow\infty.
\end{align*}
To show a similar requirement, we need to further enlarge $\mathcal{U}_{n}$. The following lemma shows how to do so.
\begin{lemma}\label{lem:LargerUn}
For any $\eta>0$, there exists $b:=b(\eta)>0$ such that for every $n\in\mathbb{N}$ large enough,
\begin{align*}
\inf_{k:\,|k-2nq|<bn}\mathbb{P}\left(U_{n}=k\right)\geq e^{-\eta n}.
\end{align*}
\end{lemma}
\noindent
\textbf{Proof.} Since $U_{n}\sim B(2n,q)$,
\begin{align*}
\mathbb{P}\left(U_{n}=k\right)=\binom{2n}{k}q^{k}(1-q)^{2n-k},\quad k=0,\ldots,2n.
\end{align*}
Next (cf.~\cite[Example 12.1.3]{CoverThomas:2006}),
\begin{align*}
\binom{2n}{k}\geq\frac{1}{2n+1}\exp\left(2n\,h_{e}\left(\frac{k}{2n}\right)\right),\quad k=0,\ldots,2n,\quad n\in\mathbb{N},
\end{align*}
where $h_{e}$ is the binary entropy function of base $e$:
\begin{align*}
h_{e}(q)=-q\ln q-(1-q)\ln(1-q),\quad q\in(0,1).
\end{align*}
Hence,
\begin{align*}
\mathbb{P}\left(U_{n}=k\right)\geq\frac{1}{2n+1}\exp\left(2n\left[h_{e}\left(\frac{k}{2n}\right)+\frac{k}{2n}\ln q+\left(1-\frac{k}{2n}\right)\ln(1-q)\right]\right).
\end{align*}
The continuous function $g(x)=h_{e}(x)+x\ln q+(1-x)\ln(1-q)$, $x\in(0,1)$, is such that $g(x)\leq 0$, for any $x\in(0,1)$, with $g(q)=0$. Thus, for any $\eta>0$, there exists $b(\eta)>0$ such that, whenever $|x-q|\leq b/2$, then $g(x)\geq -\eta/4$, and so, if $|k-2nq|\leq b(\eta)n$,
\begin{align*}
h_{e}\left(\frac{k}{2n}\right)+\frac{k}{2n}\ln q+\left(1-\frac{k}{2n}\right)\ln(1-q)>-\frac{\eta}{4}.
\end{align*}
Therefore, for $n$ (depending only on $\eta$) large enough,
\begin{align*}
\inf_{k:|k-2nq|<bn}\mathbb{P}\left(U_{n}=k\right)\geq\frac{1}{2n+1}\exp\left[2n\,h_{e}\left(\frac{k}{2n}\right)\right]q^{k}(1-q)^{2n-k}\geq\frac{1}{2n+1}\exp\left(-\frac{\eta}{2}n\right)\geq e^{-\eta n}.
\end{align*}
The proof is now complete.\hfill $\Box$

\medskip
Using Lemma \ref{lem:LargerUn}, we can now verify the condition \eqref{eq:lambda} of Theorem \ref{thm:LocUpBoundRateFunt}.
\begin{theorem}\label{thm:LowerBoundCondRateFunt}
Let $p\in(0,p_{0})$. Let $\varepsilon_{0}:=\varepsilon_{1}-\varepsilon_{2}>0$, where $\varepsilon_{1}>0$ and $\varepsilon_{2}>0$ are given in Theorem \ref{thm:LCSLargeProb}, and let $\delta:=\varepsilon_{0}/2$. Let $t_{0}>0$ be the unique positive solution to
\begin{align}\label{eq:tMq}
2\delta t-b^{2}-2\ln M_{q}(\delta t)=0.
\end{align}
Then, for any $\epsilon>0$ and any $t\in [0,t_{0}]$, there exists $N(t)\in\mathbb{N}$ such that for any $n\geq N(t)$, \eqref{eq:lambda} holds true with $\nu_{n}\equiv 1-\epsilon$, $n\in\mathbb{N}$.
\end{theorem}

\noindent
\textbf{Proof.} By Lemma \ref{lem:LargerUn}, the set $\mathcal{U}_{n}$ can now be further enlarged. Indeed, taking $0<\eta<c_{1}$ (where $c_{1}>0$ is given in Theorem \ref{thm:LCSLargeProb}), there exists $b:=b(\eta)>0$ such that
\begin{align*}
\sup_{k:\,|k-2nq|\leq bn}\frac{e^{-c_{1}n}}{\mathbb{P}\left(U_{n}=k\right)}\leq e^{(\eta-c_{1})n}\rightarrow 0,\quad n\rightarrow\infty.
\end{align*}
Let us therefore take
\begin{align*}
\mathcal{U}_{n}=\left[2nq-bn,\,2nq+bn\right]\cap\mathcal{S}_{n}.
\end{align*}
With $\varphi(n)=e^{-\eta n}$ and by Theorem \ref{thm:LCSLargeProb}, the inequalities \eqref{eq:DeltafastVarphi} (with $\varepsilon_{0}=\varepsilon_{1}-\varepsilon_{2}$) and \eqref{eq:GenLowerBoundProbU} are both satisfied. Since in our LCS case \eqref{eq:RanTransImpli} also holds, all the assumptions of Theorem \ref{thm:FastConvThm} are verified and thus, for any $k\in\mathcal{U}_{n}$, $\ell_{n}(k+1)-\ell_{n}(k)\geq\delta=\varepsilon_{0}/2$, provided that $n$ is large enough. Therefore, by the conditional Jensen's inequality,
\begin{align}
\mathbb{E}\left(e^{t\left(L_{n}(\mathbf{Z}_{n})-\ell_{n}(2nq)\right)}\right)&\geq\mathbb{E}\left(e^{t\left(\ell_{n}(U_{n})-\ell_{n}(2nq)\right)}\right)\geq\sum_{k=2nq}^{2nq+bn}e^{t\left(\ell_{n}(k)-\ell_{n}(2nq)\right)}\mathbb{P}\left(U_{n}=k\right)\nonumber\\
\label{eq:ExpLnS3Ineq} &\geq\sum_{k=2nq}^{2nq+bn}e^{t\delta\left(k-2nq\right)}\mathbb{P}\left(U_{n}=k\right)=:S_{n}^{(3)}(t).
\end{align}
On the other hand,
\begin{align*}
S_{n}(t):=M_{q}(\delta t)^{2n}&=\mathbb{E}\left(e^{\delta t\left(U_{n}-2nq\right)}\right)=\left[e^{t\delta(1-q)}q+e^{-t\delta q}(1-q)\right]^{2n}\\
&=\sum_{k=0}^{2nq-bn-1}e^{\delta t\left(k-2nq\right)}\mathbb{P}\left(U_{n}=k\right)+\sum_{k=2nq-bn}^{2nq-1}e^{\delta t\left(k-2nq\right)}\mathbb{P}\left(U_{n}=k\right)\\
&\quad+\sum_{k=2nq}^{2nq+bn}e^{\delta t\left(k-2nq\right)}\mathbb{P}\left(U_{n}=k\right)+\sum_{k=2nq+bn+1}^{2n}e^{\delta t\left(k-2nq\right)}\mathbb{P}\left(U_{n}=k\right)\\
&=:S_{n}^{(1)}(t)+S_{n}^{(2)}(t)+S_{n}^{(3)}(t)+S_{n}^{(4)}(t).
\end{align*}
It is easy to see that for every $t>0$,
\begin{align*}
S_{n}^{(2)}(t)\leq\sum_{k=2nq-bn}^{2nq-1}e^{t\delta(k-2nq)}=\sum_{j=1}^{nb}e^{-t\delta j}=\frac{e^{-t\delta}-e^{-t\delta(nb)}}{1-e^{-t\delta}}\rightarrow\frac{1}{e^{t\delta}-1}=:c(t),\quad n\rightarrow\infty.
\end{align*}
Moreover, for any $t>0$,
\begin{align*}
S_{n}^{(1)}(t)\!+\!S_{n}^{(4)}(t)\!=\!\left(\sum_{k=0}^{2nq-bn-1}\!\!\!+\!\!\!\sum_{k=2nq+bn+1}^{2n}\right)\!e^{t\delta(k-2nq)}\mathbb{P}\!\left(U_{n}\!=\!k\right)\!<\!e^{2\delta tn}\mathbb{P}\!\left(\left|U_{n}\!-\!2nq\right|\!>\!nb\right)\!\leq\!2e^{\left(2\delta t-b^{2}\right)n},
\end{align*}
where the last inequality follows from Hoeffding's exponential inequality. Thus, for $t>0$,
\begin{align*}
1=\frac{S_{n}^{(1)}(t)+S_{n}^{(2)}(t)+S_{n}^{(3)}(t)+S_{n}^{(4)}(t)}{S_{n}(t)}\leq\frac{c(t)}{S_{n}(t)}+\frac{S_{n}^{(3)}(t)}{S_{n}(t)}+2\exp\left(\left(2\delta t-b^{2}-2\ln M_{q}(\delta t)\right)n\right).
\end{align*}
Hence, for any $0<t<t_{0}$, where $t_{0}$ is the solution to \eqref{eq:tMq},
\begin{align*}
\frac{S_{n}^{(3)}(t)}{S_{n}(t)}\rightarrow 1,\quad\text{as }\,n\rightarrow\infty.
\end{align*}
Together with \eqref{eq:ExpLnS3Ineq}, it follows that, for every $t\in(0,t_{0}]$, there exists $N(t)\in\mathbb{N}$, such that for any $n\geq N(t)$,
\begin{align*}
\mathbb{E}\left(e^{t\left(L_{n}(\mathbf{Z}_{n})-\ell_{n}(2nq)\right)}\right)\geq S_{n}(t)\geq (1-\epsilon)M_{q}(\delta t)^{2n},
\end{align*}
which completes the proof.\hfill $\Box$

\section{Beyond the Binary LCS Case: Discussion}\label{sec:diss}

In this paper we applied the general methodology developed in Section \ref{sec:GenMomentBounds} to the binary LCS case with strongly asymmetric distributions. From Theorem \ref{thm:LCSLargeProb}, this setup is one of few known models where the existence of a suitable random transformation and of $\varepsilon_{0}$ (such that \eqref{eq:RanTransImpli} and $\Delta_{n}(\varepsilon_{0})\rightarrow 0$ both hold) have been proven to exist. In the binary LCS case, Theorem \ref{thm:LCSLargeProb} assumes a very asymmetric distribution, and further research is needed to relax this asymmetry assumption.

Other setups where the assumptions of Theorem \ref{thm:FastConvThm} have been proved to hold were studied in~\cite{HoudreMatzinger:2007}, \cite{LemberMatzingerTorres:2012(2)} and~\cite{LemberMatzingerTorres:2013}. For example, in~\cite{LemberMatzingerTorres:2012(2)}
and~\cite{LemberMatzingerTorres:2013}, the random variables $(X_{n})_{n\in\mathbb{N}}$ and $(Y_{n})_{n\in\mathbb{N}}$ are all i.i.d., every letter $\beta\in\mathcal{A}$ is taken with a positive probability (i.e., $\mathbb{P}(X_{1}=\beta)>0$), and the scoring function $S$ as well as the distribution of $X_{1}$ do satisfy the following {\it asymmetry assumption}.
\begin{assumption}\label{assump:mimi}
There exist $a,b\in\mathcal{A}$ such that
\begin{align}\label{eq:condmimi}
\sum_{\beta\in\mathcal{A}}\mathbb{P}(X_{1}=\beta)\left(S(b,\beta)-S(a,\beta)\right)>0.
\end{align}
\end{assumption}
For the binary alphabet $\mathcal{A}=\{a,b\}$, condition \eqref{eq:condmimi} reads
\begin{align*}
\left(S(b,a)-S(a,a)\right)\mathbb{P}\left(X_{1}=a\right)+\left(S(b,b)-S(b,a)\right)\mathbb{P}\left(X_{1}=b\right)>0.
\end{align*}
Since $S$ is symmetric and one could exchange $a$ and $b$, the condition \eqref{eq:condmimi} actually becomes
\begin{align*}
\left(S(b,a)-S(a,a)\right)\mathbb{P}\left(X_{1}=a\right)+\left(S(b,b)-S(b,a)\right)\mathbb{P}\left(X_{1}=b\right)\neq 0.
\end{align*}
When $S(b,b)=S(a,a)>S(b,a)$ (recall that $S$ is assumed to be symmetric and non-constant), then Assumption \ref{assump:mimi} is satisfied if and only if $\mathbb{P}(X_{1}=a)\neq\mathbb{P}(X_{1}=b)$. Thus, in terms of the distribution of $X_{i}$ and of the scoring function $S$, \eqref{eq:condmimi} requires much less than Theorem~\ref{thm:LCSLargeProb}. However, the price to be paid there is in terms of $\delta$. Namely, the analogue of Theorem~\ref{thm:LCSLargeProb} can be shown to hold only if the gap penalty $-\delta$ is sufficiently high. The theorem itself is as follows (cf.~\cite[Theorem 3.1]{LemberMatzingerTorres:2012(2)}):
\begin{theorem}\label{thm:app1}
Let Assumption \ref{assump:mimi} hold. Let $a,b\in\mathcal{A}$ be as in \eqref{eq:condmimi}, and let the random transformation $\mathcal{R}$ turn the letter $a$ in $\mathbf{Z}_{n}$ into the letter $b$. Then there exist constants $\delta_{0}<0$, $\epsilon_{0}>0$, $\alpha>0$, and $n_{0}\in\mathbb{N}$, such that for any $\delta<\delta_{0}$ and $n\geq n_{0}$,
\begin{align*}
\mathbb{P}\left(\mathbb{E}\left(\left.L_{n}(\mathcal{R}(\mathbf{Z}_{n}))-L_{n}(\mathbf{Z}_{n})\,\right|\mathbf{Z}_{n}\right)\geq\epsilon_{0}\right)\geq 1-e^{-\alpha n}.
\end{align*}
\end{theorem}
\begin{remark}\label{rem:GapPrice}
Typically $\delta_{0}<0$, so that the condition $\delta<\delta_{0}$ indicates that the gap penalty $-\delta$ has to be sufficiently large. Hence, the result does not apply to the LCS case. Intuitively, the larger the gap penalty (the smaller the gap price), the fewer gaps appear in the optimal alignment, so that the optimal alignment is closer to the pairwise comparison (Hamming distance). Some methods for determining a sufficient $\delta_{0}$, as well as some examples, are discussed in~\cite{LemberMatzingerTorres:2013}. We believe that the assumption on $\delta$ can be relaxed so that Theorem \ref{thm:app1} holds under more general assumptions.
\end{remark}
If the assumptions of Theorem \ref{thm:app1} hold (i.e., the scoring function $S$ and the law of $X_{1}$ satisfy Assumption \ref{assump:mimi} and $\delta$ is sufficiently small) and the alphabet $\mathcal{A}=\{a,b\}$ is binary, then the situation is exactly as in the binary LCS case considered in Section \ref{sec:LCS} except that $A$ (see \eqref{eq:DiffLtildeZLZ}) might not be $1$. However, the constant $A$ is not essential in any of the proofs, so all the results of Sections \ref{sec:LCS} and \ref{sec:Rate} continue to hold (with appropriate changes stemming from $A$). If $|\mathcal{A}|>2$, then to apply Theorem \ref{thm:app1}, one should use a more general version of Theorem \ref{thm:FastConvThm}, see~\cite[Theorem 2.2]{LemberMatzingerTorres:2012(2)}. The same theorem should apply to the LCS case when the alphabet is not binary, but the model is still strongly asymmetric, i.e., one letter has to have a very large probability.


\begin{thebibliography}{999}

\bibitem{Alexander:1994}
K. S. Alexander,
\newblock {The Rate of Convergence of the Mean Length of the Longest Common Subsequence},
\newblock {\em The Annals of Applied Probability}, Vol. 4, No. 4, pp. 1074-1082, 1994.

\bibitem{AmsaluHoudreMatzinger:2014}
S. Amsalu, C. Houdr\'{e} and H. Matzinger,
\newblock {Sparse Long Blocks and the Micro-Structure of the Longest Common Subsequences},
\newblock {\em Journal of Statistical Physics}, Vol. 154, Issue 6, pp. 1516-1549, 2014.

\bibitem{AmsaluHoudreMatzinger:2016}
S. Amsalu, C. Houdr\'{e} and H. Matzinger,
\newblock {Sparse Long Blocks and the Variance of the Length of Longest Common Subsequences in Random Words},
\newblock {arXiv:1204.1009v2}, 2016.

\bibitem{ArratiaWaterman:1994}
R. Arratia and M. S. Waterman,
\newblock {A Phase Transition for the Score in Matching Random Sequences Allowing Deletions},
\newblock {\em The Annals of Applied Probability}, Vol. 4, No. 1, pp. 200-225, 1994.

\bibitem{BonettoMatzinger:2006}
F. Bonetto and H. Matzinger,
\newblock {Fluctuations of the Longest Common Subsequence in the Asymmetric Case of 2- and 3-Letter Alphabets},
\newblock {\em Latin American Journal of Probability and Mathematical Statistics}, Vol. 2, pp. 195-216, 2006.

\bibitem{BoucheronLugosiMassart:2013}
S. Boucheron, G. Lugosi and P. Massart,
\newblock {\em Concentration Inequalities: A Nonasymptotic Theory of Independence},
\newblock {Oxford University Press}, Oxford, United Kingdom, 2013.

\bibitem{ChristianiniHahn:2007}
N. Christianini and M. W. Hahn,
\newblock {\em Introduction to Computational Genomics: A Case Studies Approach},
\newblock {Cambridge University Press}, Cambridge, United Kingdom, 2007.

\bibitem{ChvatalSankoff:1975}
V. Chv\'{a}tal and D. Sankoff,
\newblock {Longest Common Subsequences of Two Random Sequences},
\newblock {\em Journal of Applied Probability}, Vol. 12, No. 2, pp. 306-315, 1975.

\bibitem{CoverThomas:2006}
T. Cover and J. Thomas,
\newblock {\em Elements of Information Theory},
\newblock {Second Edition, John Wiley \& Sons, Inc.}, Hoboken, New Jersey, 2006.

\bibitem{DurbinEddyKroghMitchison:1998}
R. Durbin, S. Eddy, A. Krogh and G. Mitchison,
\newblock {\em Biological Sequence Analysis: Probabilistic Models of Proteins and Nucleic Acids},
\newblock {Cambridge University Press}, Cambridge, United Kingdom, 1998.

\bibitem{DurringerLemberMatzinger:2007}
C. Durringer, J. Lember and H. Matzinger,
\newblock {Deviation from the Mean in Sequence Comparison with a Periodic Sequence},
\newblock {\em Latin American Journal of Probability and Mathematical Statistics}, Vol. 3, pp 1-29, 2007.

\bibitem{GongHoudreIslak:2016}
R. Gong, C. Houdr\'{e} and \"{U}. I\c{s}lak,
\newblock {A Central Limit Theorem for the Optimal Alignments Score in Multiple Random Words},
\newblock {arXiv:1512.05699v2}, 2016.

\bibitem{GrossmannYakir:2004}
S. Grossmann and B. Yakir,
\newblock {Large Deviations for Global Maxima of Independent Superadditive Processes with Negative Drift and an Application to Optimal Sequence Alignments},
\newblock {\em Bernoulli}, Vol. 10, No. 5, pp. 829-845, 2004.

\bibitem{Hammersley:1974}
J. M. Hammersley,
\newblock {Postulates for Subadditive Processes},
\newblock {\em The Annals of Probability}, Vol. 2, No. 4, pp. 652-680, 1974.

\bibitem{HoudreIslak:2015}
C. Houdr\'{e} and \"{U}. I\c{s}lak,
\newblock {A Central Limit Theorem for the Length of the Longest Common Subsequences in Random Words},
\newblock {arXiv:1408.1559v3}, 2015.

\bibitem{HoudreMa:2016}
C. Houdr\'{e} and J. Ma,
\newblock {On the Order of the Central Moments of the Length of the Longest Common Subsequences in Random Words},
\newblock {To appear in {\em High Dimensional Probability VII: The Carg\`{e}se Volume}}, Birkh\"{a}user, 2016.

\bibitem{HoudreMatzinger:2007}
C. Houdr\'{e} and H. Matzinger,
\newblock {On the Variance of the Optimal Alignments Score for Binary Random Words and an Asymmetric Scoring Function},
\newblock {\em Journal of Statistical Physics}, Vol. 164, Issue 3, pp. 693-734, 2016.

\bibitem{KeckicVasic:1971}
J. D. Ke\v{c}ki\'{c} and P. M. Vasi\'{c},
\newblock {Some Inequalities for the Gamma Function},
\newblock {\em Publications de L'institut Math\'{e}matique}, Nouvelle S\'{e}rie, tome 11, Vol. 25, pp. 107-114, 1971.

\bibitem{LemberMatzinger:2009}
J. Lember and H. Matzinger,
\newblock {Standard Deviation of the Longest Common Subsequence},
\newblock {\em The Annals of Probability}, Vol. 37, No. 3, pp. 1192-1235, 2009.

\bibitem{LemberMatzingerTorres:2012}
J. Lember, H. Matzinger and F. Torres,
\newblock {The Rate of the Convergence of the Mean Score in Random Sequence Comparison},
\newblock {\em The Annals of Applied Probability}, Vol. 22, No. 3, pp. 1046-1058, 2012.

\bibitem{LemberMatzingerTorres:2012(2)}
J. Lember, H. Matzinger and F. Torres,
\newblock {General Approach to the Fluctuations Problem in Random Sequence Comparison},
\newblock {arXiv:1211.5072v1}, 2012.

\bibitem{LemberMatzingerTorres:2013}
J. Lember, H. Matzinger and F. Torres,
\newblock {Proportion of Gaps and Fluctuations of the Optimal Score in Random Sequence Comparison}
\newblock {In {\em Limit Theorems in Probability, Statistics and Number Theory (In Honor of Friedrich G\"otze)}}, Springer, Vol. 42, pp. 207-234,  2013.

\bibitem{LinOch:2004}
C. Y. Lin and F. J. Och,
\newblock {Automatic Evaluation of Machine Translation Quality Using Longest Common Subsequence and Skip-Bigram Statistics},
\newblock {In {\em ACL'04: Proceedings of the 42nd Annual Meeting on Association for Computational Linguistics}}, pp. 605, 2004.

\bibitem{Melamed:1995}
I. D. Melamed,
\newblock {Automatic Evaluation and Uniform Filter Cascades for Inducing N-Best Translation Lexicons},
\newblock {\em In Proceedings of the Third Workshop on Very Large Corpora}, 1995.

\bibitem{Melamed:1999}
I. D. Melamed,
\newblock {Bitext Maps and Alignment via Pattern Recognition},
\newblock {\em Computational Linguistics}, Vol. 25, Issue 1, pp. 107-130, 1999.

\bibitem{Pevzner:2000}
P. A. Pevzner,
\newblock {\em Computational Molecular Biology: An Algorithmic Approach},
\newblock {MIT Press}, Cambridge, 2000.

\bibitem{Shiryaev:1995}
A. N. Shiryaev,
\newblock {\em Probability},
\newblock {Second Edition, Springer-Verlag}, New York, 1995.

\bibitem{SmithWaterman:1981}
T. F. Smith and M. S. Waterman,
\newblock {Identification of Common Molecular Subsequences},
\newblock {\em Journal of Molecular Biology}, Vol. 147, Issue 1, pp. 195-197, 1981.

\bibitem{Steele:1986}
J. M. Steele,
\newblock {An Efron-Stein Inequality for Nonsymmetric Statistics},
\newblock {\em The Annals of Statistics}, Vol. 14, No. 2, pp. 753-758, 1986.

\bibitem{Torres:2009}
F. Torres,
\newblock {On the Probabilistic Longest Common Subsequence Problem for Sequences of Independent Blocks},
\newblock {\em Ph.D Thesis}, Bielefeld University, 2009.

\bibitem{Waterman:1994}
M. S. Waterman,
\newblock {Estimating Statistical Significance of Sequence Alignments},
\newblock {\em Philosophical Transactions of the Royal Society: Biological Sciences}, Vol. 344, Issue 1310, pp. 383-390, 1994.

\bibitem{Waterman:1995}
M. S. Waterman,
\newblock {\em Introduction to Computational Biology: Maps, Sequences and Genomes},
\newblock {Chapman \& Hall / CRC Press}, 1995.

\bibitem{YangLi:2003}
C. C. Yang and K. W. Li,
\newblock {Automatic Construction of English/Chinese Parallel Corpora},
\newblock {\em Journal of the American Society for Information Science and Technology}, Vol. 54, Issue 8, pp. 730-742, 2003.

\end{thebibliography}
\end{document}